\documentclass[twoside]{article}
\usepackage[aop]{imsart}
\usepackage{amsmath,amsthm,amssymb}
\usepackage{graphicx}
\usepackage{natbib}

\begin{document}

\renewcommand{\Re}{\mathbb{R}}
\renewcommand{\flushbottom}{false}

\begin{frontmatter}
\title{Processes with Inert Drift}
\runtitle{Processes with Inert Drift}
\author{David White}
\runauthor{David White}
\address{David White \\ Department of Mathematics \\ Cornell University \\ Ithaca, NY 14850\\ E-mail: white@u.washington.edu}

\begin{abstract}
We construct a stochastic process whose drift is a function of the process's local time at a reflecting barrier.  The process arose as a model of the interactions of a Brownian particle and an inert particle in 
\citep{knight:01}.  Interesting asymptotic results are obtained for two different arrangements of inert particles and Brownian particles.  A version of the process in $\Re^d$ is also constructed.
\end{abstract}


\end{frontmatter}

\newtheorem{thm}{Theorem}[section]
\newtheorem{lem}[thm]{Lemma}
\newtheorem{cor}[thm]{Corollary}
\theoremstyle{definition}
\newtheorem*{nota}{Notation}
\newtheorem{rem}{Remark}



\section{Introduction}
In this paper, we construct solutions $(X,L)$ to the set of equations
\begin{align*}
X(t)&=B(t)+L(t)+\int_0^t \mu(L(s))\ ds \\
X(t)&\geq 0 \\
L(t)&=\lim_{\epsilon\rightarrow 0}\frac{1}{2\epsilon}\int_0^t 1_{[0,\epsilon)}X(s)\ ds
\end{align*}
for certain monotone functions $\mu(l)$, where $B(t)$ is a standard Brownian motion with B(0)=0.  The case when $\mu(l)=Kl$, $K>0$ was studied by Frank Knight in \cite{knight:01}, which was the starting point for this research.

The case when $\mu(l)=0$ is the classical reflected Brownian motion.  There are two approaches to this case.  The first is to let $X(t)=|B(t)|$, define $L(t)$ from $X(t)$ as above, and define a new Brownian motion $\tilde B(t)=X(t)-L(t)$.  For the second approach we define
\[ L(t)= -\inf_{s<t} B(s), \]
and take $X(t)=B(t)+L(t)$.  The two approaches yield processes $(X, L)$ with the same distribution.  The second approach has the advantages of being easier to work with, and of retaining the original $B(t)$.

For other $\mu(l),$ $X(t)$ can be broken into $B(t)+L(t)$ and $-\int_0^t\mu(L(s))\ ds$. This can be interpreted as the path of a Brownian particle which reflects from the path
$-\int_0^t\mu(L(s))\ ds$ of an inert particle.  We call the second path inert because its derivative is constant when the two particles are apart.  With this model in mind, we consider two other configurations of Brownian particles and inert particles.  We also consider a generalization to $\Re^d$.

In Section 2, we consider the one-dimensional case.  We also obtain
an explicit formulas for the distribution of $\lim_{t\rightarrow\infty}\mu(L(t))$  and for the excursion measure of the process with drift,
which we use in later sections.

In Section 3, we consider the case when a 
Brownian particle and an inert particles are confined to a one-dimensional
ring.  In other words, we have a process with drift confined to a closed 
interval, with reflection at both endpoints.  
The velocity of the inert particle turns
out to be very interesting process.  Indexed by local time, it is a piecewise
linear process with Gaussian stationary distribution, though it is not a 
Gaussian process.  Some results from this section have appeared in \cite{bzwh}.
We show that under rescaling, the velocity process converges in distribution
to an Ornstein-Uhlenbeck process.

Section 4 discusses the configuration consisting of two independent Brownian
particles with an inert particle separating them.  We would like to determine
whether the process is recurrent or transient; that is, whether all three
particles can meet in finite time, or whether the inert particle in the
middle forces the distance between the Brownian particles to tend to infinity.
Interestingly, this configuration seems to be a critical case between the
two behaviors.  We show that under rescaling the distance between the two
Brownian particles is a two-dimensional Bessel process.  
Dimension
two is exactly the critical case for the family of Bessel processes.  

Section 5 extends the results in Section 2 to domains in $\Re^d$,
$d\geq 2$.  This sections uses results of Lions and Sznitman
on the existence of reflected Brownian motion in domains in $\Re^d$.
We show existence and uniqueness for the process with drift for the
case when the velocity gained is proportional to the local time.

This paper is partly based on a Ph.D.~thesis written at the University of Washington under the guidance of Chris Burdzy.

\section{Existence and Uniqueness in One Dimension}\label{ch:sk}



In this section, we prove the following version of the well-known Skorohod 
Lemma:

\begin{thm}\label{thm:esk}
Let $f(t)$ be a continuous function with $f(0)\geq 0$, and let $\mu(l)$
be a continuous monotone function with
\begin{equation}\label{eq:mucon}
 \lambda(l)=\sup_{a<b<l}
   ~\frac{|\mu(b)-\mu(a)|}{b-a}<\infty
\end{equation}
for every $l$, and if $\mu(l)\rightarrow -\infty$, then
\begin{equation}\label{eq:mucon2}
\sum_n(|\mu(n)|\vee 1)^{-1}=\infty.
\end{equation}
Then there is a unique continuous $L(t)$ satisfying
\begin{enumerate}
\newcommand{\uok}{\mu\!\circ\!L}
\item $x(t)=f(t)+L(t)+\int_0^t \uok(s)\ ds\geq 0$,
\item $L(0)=0$, $L(t)$ is nondecreasing,
\item $L(t)$ is flat off $\{t: x(t)=0\}$;
 i.e., $\int_0^\infty 1_{\{x(s)>0\}}dL(s)=0$.
\end{enumerate}
\end{thm}

As a reminder to the reader, we quote the Skorohod Lemma.  A proof can
be found in, for example, \cite{ks:91}.

\begin{thm}[The Skorohod Lemma]
Given a continuous function $f(t)$ with $f(0)\geq 0$, there
is a unique continuous function
$L(t)$ satisfying
\begin{enumerate}
\item $x(t)\equiv f(t)+L(t)\geq 0$,
\item $L(0)=0$, $L(t)$ is nondecreasing, and
\item $L(t)$ is flat off $\{t: x(t)=0\}$;
 i.e., $\int_0^\infty 1_{\{x(s)>0\}}dL(s)=0$.
\end{enumerate}
This function is given by
\begin{equation}
L(t)=\max\left[0, \max_{0\leq s\leq t}(-f(s))\right]  \label{eq:sk}
\end{equation}
\end{thm}

We will denote the
unique $L(t)$ corresponding to $f(t)$ in equation (\ref{eq:sk}) by $Lf(t)$.

A few corollaries of the Skorohod equation will be used in the proof of
Theorem \ref{thm:esk}.  We begin with these.

\begin{lem}
\label{lem:mono}
If $f(t)$ and $g(t)$ are two continuous functions with $f(0)\geq 0$, and
$f(t)\leq g(t)$ for $t\in[0,\infty)$, then $Lf(t)\geq Lg(t),\ 0\leq t<\infty$.
In particular, if $\phi(t)\geq 0$ then $L(f+\phi)(t)\leq Lf(t)$.
\end{lem}
\begin{proof}
This comes directly from equation (\ref{eq:sk}).
\end{proof}

\begin{lem}\label{lem:lfbound}
For continuous $f(t)$, $f(0)\geq0$, and $S<T$, we have
\begin{equation*}
\label{eq:lfbound}
Lf(T)-Lf(S)\leq \max\limits_{S\leq r\leq T}[f(S)-f(r)].
\end{equation*}
\end{lem}
\begin{proof}
Using equation (\ref{eq:sk}), we have
\begin{align*}
Lf(T) &= \max_{0\leq r\leq T}(-f(r)) \\
& \leq \max_{0\leq r\leq S}(-f(r)) + \max_{S\leq r\leq T}(f(S)-f(r)) \\
& \leq Lf(S)+\max_{S\leq r\leq T}[f(S)-f(r)],
\end{align*}
and the claim follows.
\end{proof}
                                                                                     
\begin{lem}\label{lem:meas}
If $f(t)$ and $g(t)$ are two continuous functions with $f(0)\geq 0$, and
$f(t)=g(t)$ for $0\leq t\leq T$, then $Lf(t)=Lg(t)$ for $0\leq t\leq T$.
\end{lem}
\begin{proof}
This also comes directly from \eqref{eq:sk}.
\end{proof}

Now we will prove Theorem \ref{thm:esk}.

\newcommand{\uok}{\mu\!\circ\!L}
\newcommand{\kt}{\tilde{L}}
\newcommand{\fe}{f^}
\newcommand{\ke}{L^\varepsilon}
\newcommand{\ken}{L^{\varepsilon(n)}}
\newcommand{\Ie}{I^\varepsilon}
\newcommand{\Ien}{I^{\varepsilon_n}}
\newcommand{\Te}{T^\varepsilon}
\newcommand{\xt}{\tilde{x}}
\newcommand{\uokt}{\mu\!\circ\!\kt}

\begin{proof}[Proof of Theorem \ref{thm:esk}]

Uniqueness will be proved first.  
Assume that both $L(t)$ and $\kt(t)$ satisfy conditions 1--3 of
Theorem \ref{thm:esk}.  
Let $Q=\inf\{t>0: L(t)\neq\kt(t)\}$,
and suppose that $Q<\infty$.  Define the function
$M(t)=\lambda(L(t)\vee\kt(t))(t-Q)$.  As the product of increasing
functions, $M(t)$ is increasing, and $M(Q)=0$.  Choose $R>Q$ so that
$M(R)<1$.
                                                                                     
By the continuity of $L(t)$ and $\kt(t)$, there exist $T$ with
$Q<T<R$, and $\delta>0$, with the property that
$|L(T)-\kt(T)|=\delta$, while $|L(t)-\kt(t)|<\delta$ for $t<T$.  We
will show that the case $L(T)-\kt(T)=\delta$ yields a contradiction;
the same argument works for the other case, $\kt(T)-L(T)=\delta$.
                                                                                     
Suppose that $L(T)-\kt(T)=\delta$ and $L(t)-\kt(t)<\delta$ for $t<T$.
We have that
\begin{align*}
 x(T)-\xt(T) &= (x(T)-x(Q))-(\xt(T)-\xt(Q)) \\
  &= [L(T)-L(Q)+\int_Q^T\uok(r)dr] \\
   &\qquad -[\kt(T)-\kt(Q)+\int_Q^T\uokt(r)dr] \\
  &= (L(T)-\kt(T))+\int_Q^T(\uok(r)-\uokt(r))dr,
\end{align*}
and that
\begin{align*}
 \left|\int_Q^T(\uok(r)-\uokt(r))dr\right|
   &\leq (T-Q)\lambda(L(T)\vee\kt(T))\sup_{Q\leq t\leq T}(L(t)-\kt(t)) \\
  &= M(T)(L(T)-\kt(T)).
\end{align*}
Combining these two equations we see that
\begin{align*}
 x(T)-\xt(T) &\geq (1-M(T))(L(T)-\kt(T)) > 0.
\end{align*}
Let $S=\sup\{t<T: x(t)=\xt(t)\}$.  Then for $S<t\leq T$, $x(t)>\xt(t)\geq 0$,
and by condition 3 of the theorem, $L(t)=L(S)$, and in particular,
$L(T)=L(S)$.  Then for $S<t\leq T$, $L(t)-\kt(t)=L(T)-\kt(t)$.  Since
$\kt(t)$ is nondecreasing by condition 2 of the theorem, $L(T)-\kt(t)$ is
nonincreasing in $t$, and we see that $L(S)-\kt(S)\geq L(T)-\kt(T)=\delta$.
This contradicts how $T$ was chosen.  Therefore, $L(t)=\kt(t)$ for all
$t\geq 0$.

Next we will construct the $L(t)$ corresponding to $f(t)$ and $\mu(t)$.
First, we may assume that $\mu(t)\geq 0$.  Otherwise, set
$m=\inf\{\mu(t), t\geq 0\}$,
$\tilde f(t)=f(t)+mt$ and $\tilde\mu(t)=\mu(t)-m$, assuming for the
moment that $m>-\infty$.
It is easily checked that the $L(t)$ satisfying the theorem for
$\tilde f(t)$ and $\tilde\mu(t)$ also
works for $f(t)$ and $\mu(t)$.

For any $\varepsilon>0$ we construct $\ke(t)$ and $\Ie(t)$ so that
$\ke(t)=L(f+\Ie)(t)$, and $\frac{d}{dt}\Ie(t)=\mu(n\varepsilon),\ 
n\varepsilon\leq \ke(t)\leq(n+1)\varepsilon$.  The construction is
recursive.  Let $\Ie_0=0$, $\ke_0(t)=Lf(t)$, and $\Te_0=\inf\{t>0:
\ke_0(t)=\varepsilon\}$.  Now for each positive integer $n$, define
\begin{align}
\Ie_{n+1}(t) &=
 \begin{cases}\label{eq:epsilons}
   \Ie_n(t), & t<\Te_n, \\
   \Ie_n(\Te_n)+\mu(n\varepsilon)(t-\Te_n), &
                \Te_n\leq t\leq\Te_{n+1},  \\
   \Ie_{n+1}(\Te_{n+1}), & \Te_{n+1}<t,
 \end{cases} \\
\ke_{n+1}(t) &= L(f+\sum_{j=1}^{n+1} \Ie_n)(t), \notag \\
\Te_{n+1}&= \inf\{t>0: \ke_n(t)=(n+1)\varepsilon\}. \notag
\end{align}
For $m\geq n$, we have that $\Ie_m(t)\geq\Ie_n(t)$, so by Lemma \ref{lem:mono},
$\ke_m(t)\leq\ke_n(t)$,
and from equation (\ref{eq:sk}), $\ke_m(t)=\ke_n(t)$
for $t<\Te_{m\wedge n}$.  Let $\ke(t)=\lim_{n\rightarrow\infty}\ke_n(t)$,
and $\Ie(t)=\lim_{n\rightarrow\infty}\Ie_n(t)$.
                                                                                     
By Lemma \ref{lem:lfbound}, we see that
\begin{align*}
\ke(t)-\ke(s) &= L(f+\Ie)(t)-L(f+\Ie)(s) \\
 & \leq \max_{s\leq r\leq t}\{(f+\Ie)(s)-(f+\Ie)(r)\} \\
 & \leq \max_{s\leq r\leq t}\{f(s)-f(r)\},
\end{align*}
that is, the family $\{\ke(t)\}_{\varepsilon>0}$ is equicontinuous.
Since $\Ie(t)\geq 0$, by Lemma~\ref{lem:mono}
$\{\ke(t)\}_{\varepsilon>0}$ is bounded for each $t$.  Ascoli--Arzel\`a
\cite{royden} then applies.
For any sequence
$\varepsilon(n)\rightarrow 0$, $\ken(t)\rightarrow L(t)$,
uniformly on compact subsets of $[0,\infty)$.
We will check that this $L(t)$ satisfies the conditions of the theorem.
That $x(t)\geq 0$ follows from the definition of $\ke_n(t)$ and
equation (\ref{eq:sk}), as does the second condition.
That $\Ien(t)\rightarrow\int_0^t \uok(s)\ ds$ follows from uniform convergence
and the inequality
\begin{align*}
\int_0^t \mu\!\circ\!\ke(t)\ ds-\Ie(t)\leq\varepsilon\lambda(t).
\end{align*}
For condition 3, notice that
$\{t: x(t)>\delta\}\subset\{t: x^{\varepsilon(n)}(t)>\delta/2\}$
for large enough $n$.
The uniqueness of $L(t)$ proved above shows that $L(t)$ is independent of the
$\varepsilon(n)$ chosen.
                                                                                     
For the case where $m=\inf\{\mu(l), l\geq 0\}=-\infty$, repeat the above 
construction
with $\mu^j(s)=\mu(s)\vee(-j)$, and denote the results as 
$L^j(t)$ and $x^j(t)$.
By lemma \ref{lem:meas}, $L^j(t)$ will agree with $L(t)$ for $0\leq T^j$, where
$T^j=\inf\{t: \uok(t)=-j\}$.  To complete the proof, it is only necessary 
to show that
$T^j\rightarrow\infty$.  Suppose that $T^j\rightarrow T$.  There are then 
$T_n\uparrow T$
so that $L(T_n)=n$ and $x(T_n)=0$.  We compute that
\begin{align*}
f(T_n)-f(T_{n+1}) &= 1 + \int_{T_n}^{T_{n+1}}\uok(s)\ ds.
\end{align*}
By the continuity of $f(t)$, the LHS approaches $0$, so that
\begin{align*}
\left|\int_{T_n}^{T_{n+1}}\uok(s)\ ds\right|\rightarrow 1.
\end{align*}
However,
\begin{align*}
 \left|\int_{T_n}^{T_{n+1}}\uok(s)\ ds\right| &\leq
   (T_{n+1}-T_n)\sup_{n\leq l\leq n+1}|\mu(l)|, \\
 \intertext{so that}
 T_{n+1}-T_n &\geq \frac{0.9}{|\mu(n+1)|}
\end{align*}
for sufficiently large $n$.  Then $T=T_1+\sum(T_{n+1}-T_n)\geq\infty$ by
\eqref{eq:mucon2}, a contradiction.  Hence $T^j\rightarrow\infty$, and the
proof is complete.
\end{proof}

The following corollary restates Theorem \ref{thm:esk} to be similar to Knight's original result.
Theorem \ref{thm:esk} yields a process with drift reflected at $0$, this version yields a process reflected from a second curve $v(t)$.

\begin{cor}\label{cor:esk}
Let $f(t)$ be a continuous function with $f(0)\geq 0$, and let $\mu(l)$
be a continuous increasing function with
\[ \lambda(l)=\sup_{a<b<l}
   ~\frac{|\mu(b)-\mu(a)|}{b-a}<\infty\]
for every $l\geq 0$.  Then there is a unique continuous $L(t)$ satisfying
\begin{enumerate}\label{cor:esk}
 \item $x(t)=f(t)+L(t)$,
 \item $L(0)=0$, $L(t)$ is nondecreasing,
 \item $L(t)$ is flat off $\{t: x(t)=v(t)\}$, and
 \item $v(0)=0$, $\frac{d}{dt}v(t)=-\uok(t)$.
\end{enumerate}
\end{cor}
\begin{proof}
This is just a restatement of Theorem \ref{thm:esk}.  To see this, let
\begin{align*}
v(t) &= -\int_0^t\uok(s)\ ds.
\end{align*}
\end{proof}

\begin{figure}[h]
\centering
\includegraphics[width=4in]{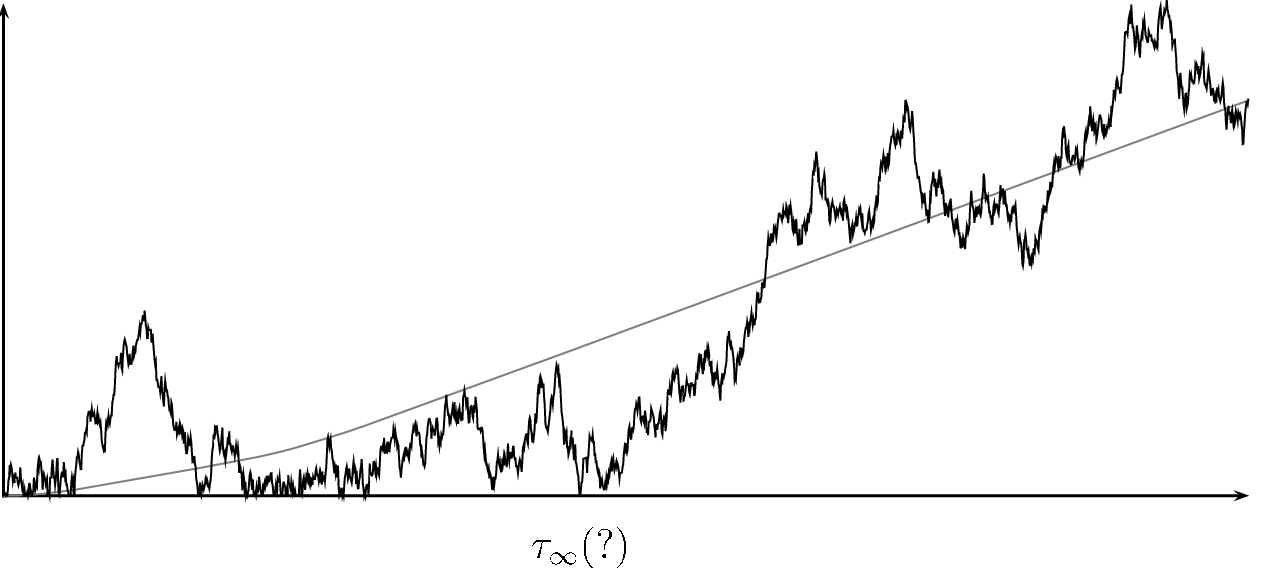}
\includegraphics[width=4in]{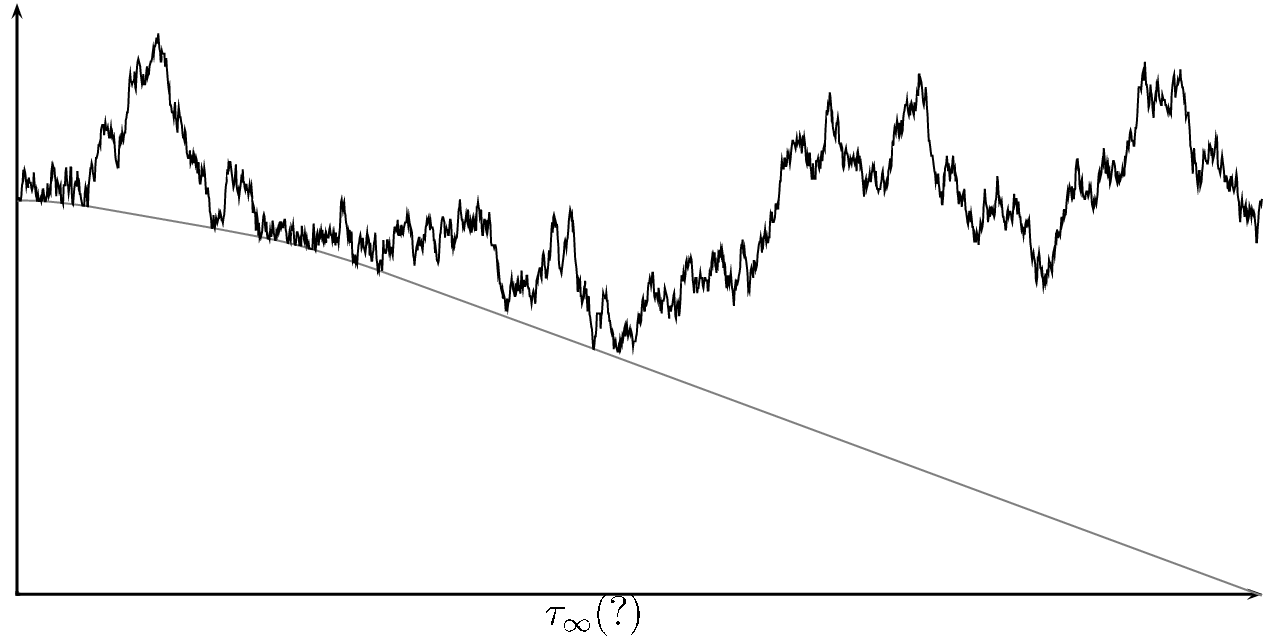}
\caption[Two equivalent versions of the process]{Two equivalent versions
of the process, corresponding 
to Theorem \ref{thm:esk} and Corollary \ref{cor:esk}, resp. }
\end{figure}

The remarks that follow demonstrate the necessity of the restrictions
(\ref{eq:mucon}) and (\ref{eq:mucon2})
on $\mu(t)$ in Theorem \ref{thm:esk}.

                                                                                     
\begin{rem}[Non-uniqueness of $L(t)$]\label{rem:ex}
In the context of Theorem~\ref{thm:esk}, if we let
$ f(t)=-t$, $\mu(l) = 1-\sqrt{l}$, $L(t)= 0$, and
$\kt(t) = t^2/4$, then $x(t)=\tilde x(t)=0$, which shows that
the theorem cannot be extended to general $f(t)$ and $\mu(l)$.
\end{rem}
                                                                                     
\begin{rem}[Blow--up of $L(t)$]
Again, in the context of Theorem~\ref{thm:esk}, we let
$f(t)=-t$, $\mu(l)=-l^2$, and $L(t)=\tan t$.  Then
\begin{align*}
x(t) &= -t+\tan t-\int_0^t\tan^2 s\ ds = 0,
\end{align*}
so $L(t)$ satisfies the conditions of the theorem, but blows up at $t=\pi/2$.
\end{rem}


\newcommand{\indic}{\mathbf 1}
\newcommand{\as}{a.s.}

The results so far do not rely on probability.  We will now remedy that by applying the results to Brownian paths.
For the rest of the section, we will need a standard Brownian motion
$\{B_t, \mathcal F_t\}$, with $B_0=0$ a.s., and some fixed $\mu(t)$
satisfying (\ref{eq:mucon}).  In the statements below, $L(t)$ will
really be $L(\omega,t)=LB(\omega,t)$.

\begin{thm}\label{thm:ltime}
If $f(t)$ in Theorem \ref{thm:esk} is replaced with a Brownian motion $B_t$,
then the
corresponding $L(t)$ is the semimartingale local time at zero, $\Lambda_t(0)$,
 \cite[p.~218]{ks:91}
of $X_t\equiv B_t+L(t)+\int_0^t\uok(s)ds$, \as

Further, we get that, \as,
\begin{align*}
2L(t) &= \lim_{\varepsilon\rightarrow 0} \varepsilon^{-1}
  \int_0^t \indic_{[0,\varepsilon)}\left(B_s+L(s)+\int_0^s\uok(r)dr\right)ds.
\end{align*}
\end{thm}
\begin{proof}
To see that $L(t)=\Lambda_t(0)$, use a version of Ito's formula
(e.g.~(7.4) on p.~218 of \cite{ks:91}),
with the identifications $M_t=B_t$, $V_t=L(t)+\int_0^t\uok(s)ds$, and
$f(x)=|x|$, and get that 
\begin{align*}
  |X_t| &= X_0-\int_0^t\indic_{\{0\}}(X_s)dB_s+
  \int_0^t\indic_{(0,\infty)}(X_s)dB_s-\int_0^t\indic_{\{0\}}(X_s)dk_s
  \\ &\qquad   +\int_0^t\uok(s)ds+2\Lambda_t(0)\\
  &= B_t-L(t)+\int_0^t\uok(s)ds+2\Lambda_t(0) \\
  \intertext{which, since $X_t\geq 0$,} &= B_t+L(t)+\int_0^t\uok(s)ds.
\end{align*}
Therefore, $L(t)=\Lambda_t(0)$, \as  

Having established the previous fact, we use properties of $\Lambda_t(a)$
to prove the other assertion.  References, again, are to \cite[p.~218]{ks:91}.
By (7.3) we get
\begin{align*}
\lim_{\varepsilon\rightarrow0}\varepsilon^{-1}\int_0^t\indic_{[0,\varepsilon)}
      \left(B_s+L(s)+
      \int_0^s\uok(r)dr\right)ds &= 
       \lim_{\varepsilon\rightarrow0}\frac{2}{\varepsilon}
        \int_0^\varepsilon\Lambda_t(a)da.
\end{align*}
By (iv), we have that, a.s., $\lim\limits_
 {a\rightarrow0}\Lambda_t(a)\rightarrow\Lambda_t(0)$, so the above limit is
$2\Lambda_t(0)$, or $2L(t)$.  
\end{proof}

\begin{cor}\label{cor:levy}
For the case $\mu\equiv\lambda$, that is, Brownian motion with drift 
$\lambda$, reflected at $0$,
$X_t\equiv B_t+L_t+\lambda t$,
\begin{align*}
  L_X(t)\equiv L_t=\max_{0\leq s\leq t} (-B_t-\lambda t).
\end{align*}
\end{cor}
\begin{proof}
This follows directly from the previous theorem and equation (\ref{eq:sk})
on page \pageref{eq:sk}.  By the notation $L_X(t)$ we will mean the
semimartingale local time at $0$ of the reflected process $X(t)$.  Note that in the stochastic case, unlike the deterministic case, the $L_X(t)$ can be recovered almost surely from $X(t)$.
\end{proof}

We now fix $\lambda$, and define $X(t)\equiv B(t)+L_X(t)+\lambda t$ as above.
Let $\phi$ be the right-continuous inverse of $L_X(t)$; that is,
\begin{equation*}
\phi(b)=\inf\{t: L_X(t)=b\}.
\end{equation*}
By Corollary \ref{cor:levy}, we have that
\begin{equation*}
\phi(b)\overset{d}{=}T_b,
\end{equation*}
where $T_b=\inf\{t: -B_t-\lambda t=b\}$.  From \cite[p.~196]{ks:91}
we have that
\begin{equation*}
P[T_b-T_a\in dt]=\frac{b-a}{\sqrt{2\pi t^3}}
 \exp\left[{-\frac{(b-a)^2}{2t}}\right]dt,
\end{equation*}
and a computation shows that for $\mu\equiv\lambda$,
\begin{equation*}
\lim_{\varepsilon\rightarrow0}\frac{1}{\varepsilon}
  P[T_{b+\varepsilon}-T_b\in dt]=\frac{dt}{\sqrt{2\pi t^3}}
   \exp\left[-\frac{\lambda^2 t}{2}\right].
\end{equation*}
Additionally, we have for $b>0$ that
\begin{align}\label{eq:infty}
P\{T_b<\infty\} &=
 \begin{cases}
  1, & \lambda\geq 0, \\ \mathrm e^{-2|\lambda| b}, &\lambda < 0,
 \end{cases} \notag \\
 &= \mathrm e^{-2|\lambda\wedge 0| b}.
\end{align}

\begin{thm}\label{thm:escape}
Let $\tau_\infty=\inf\{s: \phi(s)=\infty\}$.  Then
\begin{align*}
P\{\tau_\infty>\tau\}=\exp(-2\int_0^\tau |\mu(s)\wedge 0|ds).
\end{align*}
\end{thm}
\begin{proof}
\newcommand{\ti}{\tau_\infty}
\newcommand{\tei}{\tau^\varepsilon_\infty}
\newcommand{\tu}{\underline \tau}
\newcommand{\ddt}{\frac{d}{d\tau}}
An equivalent definition of $\tau_\infty$ is 
$\ti=\sup\{L_X(t): t\geq 0\}$.  Recall from the proof of Theorem
\ref{thm:esk} the definition of $\ke(t)$, and let 
$\tei=\sup\{\ke(t): t\geq 0\}$.  It is clear from the uniform convergence
$\ke(t)$ to $L_X(t)$ that $\ti=\lim\limits_{\varepsilon\rightarrow 0}\tei$.  

For notational convenience, define 
$\tu=\sup\{j\varepsilon: j\varepsilon\leq \tau, j \text{ an integer}\}$,
where the value of $\varepsilon$ should be apparent from the context.  
From \eqref{eq:epsilons} and \eqref{eq:infty} we see that
\begin{align*}
P(\tei>\tau|\tei>\tu)=\mathrm e^{-2|\mu(\tu)\wedge 0|(\tau-\tu)}.
\end{align*}
Using this, compute
\begin{align*}
\ddt P(\tei>\tau) &= \lim_{\delta\rightarrow0}
     \delta^{-1} [P(\tei>\tau+\delta)-P(\tei>\tau)] \\
 &= \lim_{\delta\rightarrow0} 
    \delta^{-1} [P(\tei>\tau)P(\tei>\tau+\delta|\tei>\tau)-P(\tei>\tau)] \\
 &= \lim_{\delta\rightarrow0}
    P(\tei>\tau)\delta^{-1}(\mathrm e^{-2|\mu(\tu)\wedge 0|\delta}-1) \\
 &= P(\tei>\tau)(-2|\mu(\tu)\wedge 0|).
\end{align*}
Solving this separable differential equation gives 
\begin{align*}
P(\tei>\tau)=\exp(-2\int_0^\tau|\mu(\underline s)\wedge 0|\ ds),
\end{align*}
and taking the limit as $\varepsilon\rightarrow 0$ gives the result.  
\end{proof}



For calculations, it is useful to consider the process we have constructed
as a point process of excursions the origin, indexed by local time.  Here
we derive some of the formulas that will be used in later sections.

\begin{nota}
For a process $X_t$ and a subset $C$ of $[0,\infty)\times[0,\infty]$, we define
\begin{equation*}
\nu(X; C)=\#\left\{(\tau,\lambda)\in C: \phi(\tau)-\phi(\tau-)=\lambda\right\}.
\end{equation*}
We also define the measure
\begin{equation*}
n(\cdot{})=\mathrm E \nu(X; \cdot{}).
\end{equation*}
\end{nota}

\newcommand{\dt}{d\tau}
\newcommand{\nt}{n_\tau}

\begin{thm}\label{thm:exdens}
The measure $n(\cdot)$ has density function 
\begin{equation*}
  \frac{1}{\sqrt{2\pi \lambda^3}}\exp\left(\frac{-\mu^2(\tau)\lambda}{2}\right)
     \exp\left(-2\int_0^\tau|\mu(s)\wedge 0|ds\right).
\end{equation*}
\end{thm}
\begin{proof}
For a fixed $\tau$ and $\lambda$, this decomposes as the product of the
probability that an excursion of Brownian motion with drift $\mu(\tau)$
has time duration $\lambda$ given that 
$\tau_\infty>\tau$, times the probability that
$\tau_\infty>\tau$ (from Theorem \ref{thm:escape}).
If $\tau_\infty<\tau$, no more excursions occur.
\end{proof}



The final calculation in this section is probability that an excursion of Brownian motion with constant drift $\mu$ reaches a certain height.  We will need this in the sections that follow.


\begin{lem}\label{lem:lrate}
For Brownian motion with drift $\mu$,
the intensity measure of excursions that reach level $l$ 
before returning to the origin
is given by
\begin{equation*}
\lim_{\varepsilon\rightarrow 0}\frac{1}{\varepsilon}P^{(\mu)}_\varepsilon(T_l<T_0)=
\frac{\mu\mathrm e^{\mu l}}{\sinh(\mu l)}.
\end{equation*}
\end{lem}
\begin{proof}
Apply the Girsanov transform to get that
\begin{equation*}
P^{(\mu)}_\varepsilon(T_l<T_0\wedge t)=E[1_{\{T_l<T_0\wedge t\}}Z_t],
\end{equation*}
where $Z_t=\exp(\mu (B_t-B_0)-\mu^2 t/2)$.  By the Optional Sampling theorem, 
this is
\begin{align*}
&=E[1_{\{T_l<T_0\wedge t\}}Z_{T_l\wedge t}] \\
&= E[1_{\{T_l<T_0\wedge t\}}Z_{T_l}] \\
&= \exp(\mu (l-\varepsilon))E_\varepsilon[\exp(-\mu^2 T_l/2)1_{\{T_l<T_0\}}] \\
&= \exp(\mu (l-\varepsilon))\frac{\sinh(\varepsilon\mu)}{\sinh(l\mu)}.
\end{align*}
The last formula comes from \cite[2.8.28]{ks:91}.
\end{proof}



\section{A Process with Inert Drift in an Interval}


The construction method in Section \ref{ch:sk} for one inert particle and
one Brownian particle can be extended to other configurations of particles.
In this section we construct
a process $X(t)$, which is
a Brownian motion confined to the interval $[0,l]$, with drift 
$V(t)=K(L_0(t)-L_l(t))$, where $L_0(t)$ and $L_l(t)$ are the local times at 
$0$ and $l$, resp., accumulated before time $t$.

Another way to look at this process is to allow the interval to move,
that is, construct processes $Y_0(t)\leq X(t)\leq Y_l(t)$, where $Y_0(t)$
and $Y_l(t)$ describe the paths of inert particles with drift
$V(t)=K(L_0(t)-L_1(t))$, with $Y_l(t)=Y_0(t)+l$, and with $X(t)$ the
path of a Brownian particle reflected from $Y_0(t)$ and $Y_l(t)$.
If we look at these processes modulo $l$, then $Y_0(t)$ and $Y_l(t)$ can be
seen as two sides of a single inert particle, on the boundary of a disc.
If we let $l=2\pi$, then $\exp(i X(t))$ and $\exp(i Y_0(t))$ trace out
the paths of
one Brownian particle and one inert particle on the boundary of 
the unit disc.

\begin{thm}
Given constants $l,K>0$, $x\in[0,l]$, $v\in\Re$, and a Brownian motion $B(t)$ 
with $B(0)=0$, there exist unique processes $L_0(t)$ and $L_l(t)$
satisfying
\begin{enumerate}
\item $Y_0(t)\leq X(t)\leq Y_l(t)$, where
\begin{enumerate}
\item $X(t)\equiv B(t)+L_0(t)-L_l(t)$
\item $V(t)\equiv v-K(L_0(t)-L_l(t))$
\item $Y_0(t)\equiv \int_0^t V(s)\,ds$
\item $Y_l(t)\equiv Y_0(t)+l$
\end{enumerate}
\item $L_0(t)$ and $L_l(t)$ are continuous, nondecreasing functions,
with $L_0(0)=L_l(0)=0$, and 
\item $L_0(t)$ and $L_l(t)$ are flat
off the sets $\{t:Y_0(t)=X(t)\}$ and $\{t:Y_l(t)=X(t)\}$ , resp.
\end{enumerate}
\end{thm}

\begin{proof}
  In Section \ref{ch:sk} we constructed a similar process for just one
  Brownian and one inert particle.  Because the two inert particles in
  this theorem are always distance $l$ apart, we can carry out the
  construction piecewise; that is, do the construction for one inert
particle and one Brownian particle as in the previous section until
the distance between these two particles reaches $l$, then continue
the construction with the Brownian particle and the other inert particle,
until the distance between these two reaches $l$, and continue switching
between them.
  
  Using Corollary \ref{cor:esk}, with the identifications
  $\mu(\tau)=K\tau$, $f(t)=B(t)$, we first construct unique processes
  $Y_0(t)$, $X(t)$, and $L_0(t)$ having the properties stated in the
  theorem, for $0\leq t<T_0\equiv\inf\{t:X(t)-Y_0(t)=l\}$.  Applying
  Corollary \ref{cor:esk} again, but changing the order of the
  Brownian and inert paths, we can extend $Y_0(t)$, $X(t)$, and
  $L_l(t)$ to the time interval $0\leq
  t<T_1\equiv\inf\{t>T_0:Y_l(t)-X(t)=l\}$.  Repeat this process to
extend $Y_0(t)$ and $X(t)$ until $T_2\equiv\inf\{t>T_1:X(t)-Y_0(t)=l\}$.
Then we are in the initial situation of the inert particle in contact with
the lower upper inert particle $Y_l(t)$, and we can repeat the steps
of the construction.
\end{proof}
\begin{figure}[ht]
\centering
\includegraphics[width=4in]{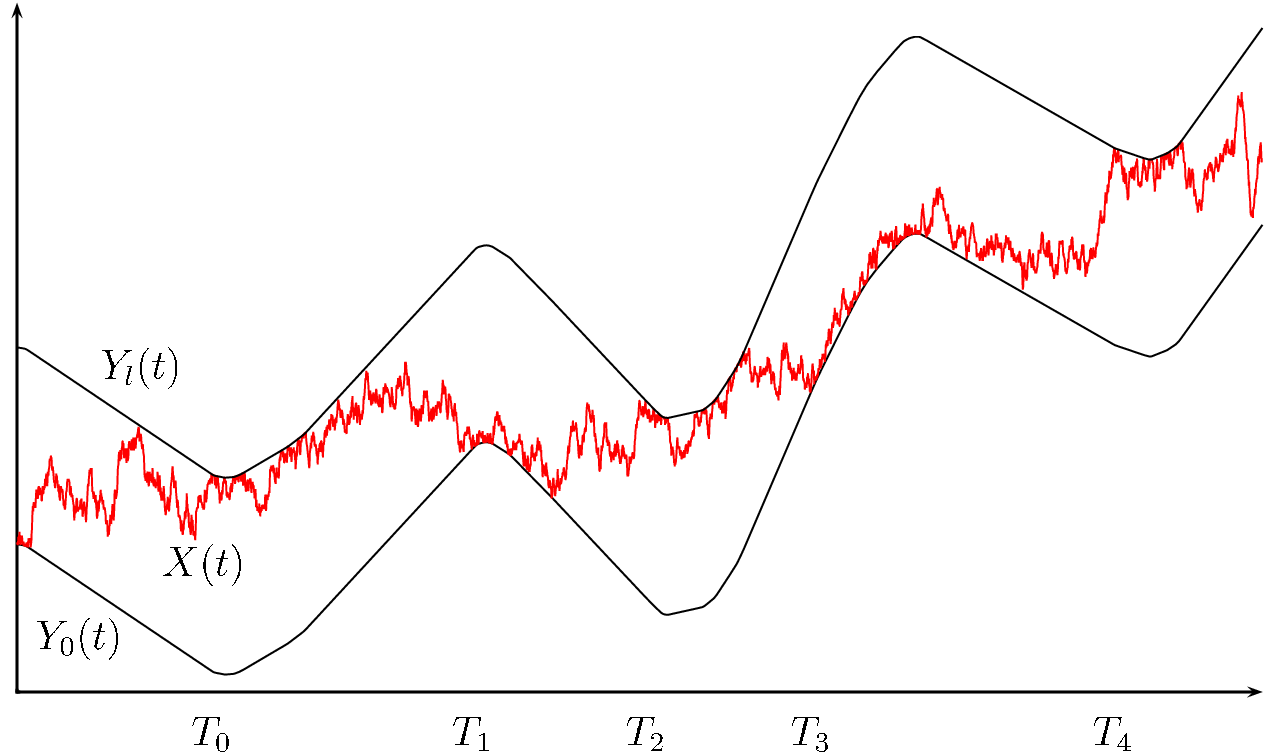}
\caption{A process with drift in an interval}
\end{figure}

\subsection{Density of the velocity process}

\newcommand{\vtau}{\tilde{V}_\tau}

The function $L(t)\equiv L_0(t)+L_l(t)$ is the total semimartingale
local time that the Brownian particle $X(t)$ spends at the endpoints
of the interval $[0,l]$.  Since the velocity process $V(t)$ changes
only at these endpoints, we will reparametrize $V$ in terms of the
local time.  We define $\vtau\equiv V(L^{-1}(\tau))$.

The process $\vtau$ is a piecewise linear process consisting of
segments with slope $\pm K$.  The slope of the process at a particular
time $\tau$ depends on which endpoint $X$ has most recently visited.
For this reason, $\vtau$ is not Markov.  We therefore introduce a
second process $J_\tau$, taking values $0$ or $1$, with value $0$
indicating that $X$ has most recently visited endpoint $0$ and the
velocity is increasing, and value $1$ indicating that $X$ has most
recently visited endpoint $l$ and the velocity is decreasing.
This technique of introducing an additional state process is similar
to that used in \cite{bkaspi}.

\newcommand{\ddv}{\frac{\partial}{\partial v}}

\begin{lem}\label{lem:vjgen}
The generator, $A$, of the process $(\vtau, J_\tau)$ is given by
\begin{align*}
Af(v,1) &= -K\frac{\partial}{\partial
v}f(v,1)+\frac{v\mathrm{e}^{-vl}}{\sinh(vl)}[f(v,0)-f(v,1)], \\
Af(v,0) &= K\frac{\partial}{\partial
v}f(v,0)+\frac{v\mathrm{e}^{vl}}{\sinh(vl)}[f(v,1)-f(v,0)]. 
\end{align*}
\end{lem}
\begin{proof}
We will prove the lemma using
the definition of the generator and Lemma \ref{lem:lrate}.
\[
Af(v,j)=\lim_{\tau\rightarrow 0}\frac{E^{(v,j)}[f(\vtau,J_\tau)-f(v,j)]}{\tau}.
\]
By Lemma \ref{lem:lrate}, excursions from $0$
reach $l$ 
with Poisson rate $v\mathrm{e}^{vl}/\sinh(vl)$, and by symmetry,
excursions from $l$ that reach $0$ occur with rate
$v\mathrm{e}^{-vl}/\sinh(vl)$.  Let $\sigma$ be the time of the first
crossing.  We can rewrite the previous equation for $j=0$ as
\begin{align*}
  Af(v,0)&=\lim_{\tau\rightarrow 0}\{\frac 1\tau P^{(v,0)}(\sigma>\tau)[f(v+K\tau,0)-f(v,0)] \\
  &\qquad+ \frac 1\tau P^{(v,0)}(\sigma<\tau)[f(v+K\sigma-K(\tau-\sigma),1)-f(v,0)]\} \\
  &= \lim_{\tau\rightarrow 0}\{P^{(v,0)}(\sigma>\tau)\frac 1\tau[f(v+K\tau,0)-f(v,0)] \\
  &\qquad+ \frac 1\tau
  P^{(v,0)}(\sigma<\tau)[f(v+K(2\sigma-\tau),1)-f(v,0)]\}.
\end{align*}
Because the limit of $P^{(v,0)}(\sigma<\tau)/\tau$ is exactly the Poisson rate given above,
$Af(v,0)$ is as stated in the lemma.  A similar calculation gives $Af(v,1)$.
\end{proof}

\begin{lem}\label{lem:adjoint}
The (formal) adjoint, $A^{*}$ of the generator $A$ is given by
\begin{align*}
A^{*}g(v,1) &= K\frac{\partial}{\partial v}g(v,1)+ 
\frac{v\mathrm{e}^{vl}}{\sinh(vl)}g(v,0)-
\frac{v\mathrm{e}^{-vl}}{\sinh(vl)}g(v,1), \\
A^{*}g(v,0) &= -K\frac{\partial}{\partial v}g(v,0)-
\frac{v\mathrm{e}^{vl}}{\sinh(vl)}g(v,0)+
\frac{v\mathrm{e}^{-vl}}{\sinh(vl)}g(v,1).
\end{align*}
\end{lem}
\begin{proof}
The (formal) adjoint $A^*$ of $A$ is the operator satisfying, for all 
suitable $f, g$
\begin{equation}
\int (Af) g dv = \int f (A^*g) dv.
\end{equation}
Let us assume that $A$ is of the somewhat more general form
\begin{align*}
Af(v,1) &= -K\frac{\partial}{\partial v}f(v,1)+a(v)[f(v,0)-f(v,1)], \\
Af(v,0) &= K\frac{\partial}{\partial v}f(v,0)+b(v)[f(v,1)-f(v,0)]. 
\end{align*}
Then we have that
\begin{align*}
\int(Af)g &= \int -K\ddv f(v,1)g(v,1)+\int a(v)[f(v,0)g(v,1)-f(v,1)g(v,1)] \\
 &+ \int K\ddv f(v,0)g(v,0)+\int b(v)[f(v,1)g(v,0)-f(v,0)g(v,0)]. \\
\intertext{Integrating by parts,}
\int(Af)g &= \int\{Kf(v,1)\ddv g(v,1)+f(v,1)[b(v)g(v,0)-a(v)g(v,1)] \} \\
 &+ \int\{-Kf(v,0)\ddv g(v,0)+f(v,0)[-b(v)g(v,0)+a(v)g(v,1)] \}. 
\end{align*}
Factoring out $f(v,j)$ leaves
\begin{align*}
A^*g(v,1) &= K\ddv g(v,1)+b(v)g(v,0)-a(v)g(v,1), \\
A^*g(v,0) &= -K\ddv g(v,0)-b(v)g(v,0)+a(v)g(v,1).
\end{align*}
\end{proof}

\begin{lem}
The process $(\vtau, J_\tau)$ has stationary normal density 
\begin{equation}
g(v,j) = \frac{1}{2\sqrt{\pi K}}\mathrm{e}^{-v^2/K}\ dv
\end{equation}
\end{lem}
\begin{proof}
The stationary distribution, $\mu$, of a process is the probability 
measure satisfying
\begin{equation}
\int Af\ d\mu=0
\end{equation}
for all $f$ in the domain of $A$.  If we assume that $d\mu$ is of the
form $g(v,j)dv$, then this is equivalent to
\begin{equation}
\int f A^{*}g\ dv=0,
\end{equation}
so that it is sufficient to find $g(v,j)$ satisfying $A^{*}g(v,j)=0$.  By
Lemma \ref{lem:adjoint}, 
\begin{align*}
A^{*}g(v,1)+A^{*}g(v,0)=K(\frac{\partial}{\partial
v}g(v,1)-\frac{\partial}{\partial v}g(v,0)), 
\end{align*}
so that $g(v,0)$ and $g(v,1)$
differ by a constant.  Note that this does not depend on the jump intensities
$a(v)$ and $b(v)$.
Since $\int(g(v,1)+g(v,0))dv =1$, $g(v,0)=g(v,1)$.
Using this fact and Lemma \ref{lem:adjoint}, we get 
\begin{align*}
A^*g(v,1) &= K\ddv g(v,1)+\frac{v(\mathrm e^{vl}-\mathrm e^{-vl})}{\sinh(vl)}g(v,1) \\
 &= K\ddv g(v,1)+2v g(v,1).
\end{align*}
This gives the separable differential equation
\[ 0=K\ddv g(v,1)+2v g(v,1), \]
which has solutions of the form $g(v,1)=C\exp(-v^2/K)$.

\end{proof}




\subsection{Behavior as $l\rightarrow 0$}

\newcommand{\vlt}{\tilde{V}^{l}(\tau)}
\newcommand{\tl}{\tau^l}
\newcommand{\sgl}{\sigma^l}
\newcommand{\dl}{{q}^l}
\newcommand{\ul}{{p}^l}
\newcommand{\dls}{\tilde{q}^l}
\newcommand{\uls}{\tilde{p}^l}
\newcommand{\vtw}{\hat{V}^l}

In this section we will let $\vlt$ denote the process $\vtau$ constructed in the 
previous section, with the constant $K=1$.  
We will show that the sawtooth process $\vlt$ converges weakly to
the Ornstein-Uhlenbeck process, under appropriate rescaling, as $l\rightarrow 0$.  The proof 
uses results from section 11.2 of \cite{sv}.

We first consider the sawtooth process $\vlt$ only at those times when it 
switches from decreasing to increasing, call them $\tl_0=0, \tl_1, \tl_2, \dots$
We will also refer to those times when the process switches from increasing
to decreasing, call them $\sgl_1, \sgl_2, \dots$  Let
$\ul_n=\sgl_n-\sgl_{n-1}$, and $\dl_n=\tl_n-\sgl_n$.  


Next, we construct the 
piecewise linear processes $\vtw(t)$ and $T^l(t)$.  We define $\vtw(t)$ so that
$\vtw(nl^2)=\tilde{V}^l(\tl_n)$, and $T^l(nl^2)=l\tl_n$.  


As in the previous section, we will let $a(v,l)=v\mathrm{e}^{vl}/\sinh(vl)$
and $b(v,l)=v\mathrm{e}^{-vl}/\sinh(vl)$.  We will use the following 
properties of $a$ and $b$ in our proof:
\begin{align*}
0 & \leq \frac{\partial}{\partial v} a(v,l) \leq 2 \\
-2 & \leq \frac{\partial}{\partial v} b(v,l) \leq 0,
\end{align*}
so that $a(v,l)$ and $b(v,l)$ are monotone functions.  We also use the
property that 
$l a(v,l)\rightarrow 1$ and $l b(v,l)\rightarrow 1$ as $l\rightarrow 0$:
\[
\lim_{l\rightarrow 0}l a(v,l)=\lim_{l\rightarrow 0}\frac{lv e^{vl}}{\sinh vl}
=\lim_{l\rightarrow 0}\frac{v(vl+1)e^{vl}}{v\cosh vl}=1.
\]



Let $P^l_{v,t}$ be the probability measure 
corresponding to $\langle\vtw(t),T^l(t)\rangle$ starting from $(v,t)$.  Let 
$P_{v,t}$ be the unique solution of the martingale problem for
\[ L= \frac{\partial^2}{\partial v^2}-2v\frac{\partial}{\partial v}
     +2\frac{\partial}{\partial T} \]
starting from $(v,t)$.
\begin{lem}\label{lem:genconv}
$P^l_{v,t}\rightarrow P_{v,t}$
as $l\rightarrow 0$ uniformly on compact subsets of $\Re^2$.
\end{lem}

The proof of Lemma \ref{lem:genconv}
depends on and will follow a number of technical lemmas.
We begin with some definitions.
The notation $P^{(v,t)}(dr,ds)$ 
denotes the transition density of $(\ul, \dl)$.
\begin{align*}
E_l &= \{(u,d)| (u-d)^2+l^2(u+d)^2\leq 1\}, \\
a_l^{VV}(v,t) &= \frac{1}{l^2}\int_{E_l}(r-s)^2 P^{(v,t)}(\ul_1\in dr,\dl_1\in ds), \\
a_l^{VT}(v,t) &= \frac{1}{l^2}\int_{E_l}l(r+s)(r-s) P^{(v,t)}(\ul_1\in dr,\dl_1\in ds), \\
a_l^{TT}(v,t) &= \frac{1}{l^2}\int_{E_l}l^2(r+s)^2 P^{(v,t)}(\ul_1\in dr,\dl_1\in ds),\\
b_l^V(v,t) &= \frac{1}{l^2}\int_{E_l}(r-s) P^{(v,t)}(\ul_1\in dr,\dl_1\in ds), \\
b_l^T(v,t) &= \frac{1}{l^2}\int_{E_l}l(r+s) P^{(v,t)}(\ul_1\in dr,\dl_1\in ds), \\
\Delta_l^{\varepsilon}(x)&=\frac{1}{l^2}P^{(v,t)}((v^l_1, l \tl_1)>\varepsilon).
\end{align*}

To prove the lemma, we need to show that the quantities
$||a_l^{VV}(v,t)-2||$, $||a_l^{VT}(v,t)||$, $||a_l^{TT}(v,t)||$,
$||b_l^V(v,t)+2v||$, $||b_l^T(v,t)-2||$ and $||\Delta_l^{\varepsilon}(v,t)||$
converge to zero as $l\rightarrow 0$, uniformly for $|v|\leq R$.  When this is
done, Lemma 11.2.1 of \cite{sv} completes the proof.

Because the density $P^{(v,t)}(dr,ds)$ is a bit unwieldy for direct computation,
we introduce $\uls$ and $\dls$,
which are exponential random variables with rate $a(v,l)$ and $b(v,l)$,
respectively.  Define
$\tilde{a}_l^{VV}(v,t)\dots\tilde{b}_l^T(v,t)$ for $\uls$ and $\dls$ as above,
using $\tilde P^{(v,t)}(dr,ds)$ for the transition density of $(\uls, \dls)$.

\begin{rem}
We will need to compare $P^{(v,t)}(dr,ds)$ and $\tilde P^{(v,t)}(dr,ds)$ in
the calculations that follow.
If the process $\vtw(t)$ stays between $v-\Delta$ and $v+\Delta$ for 
$0\leq t\leq l^2$, then the rate of the process $\ul$ is bounded by 
$a(v,l)$ and $a(v,l)+2\Delta$, and the rate of the process $\dl$ is bounded by
$b(v,l)-2\Delta$ and $b(v,l)+2\Delta$.  Then the densities satisfy
\begin{multline*}
a(v,l)(b(v,l)-2\Delta))e^{(a(v,l)+2\Delta)r}e^{(b(v,l)+2\Delta)s}
 \\ \leq  P^{(v,t)}(dr,ds)  \leq
(a(v,l)+2\Delta)(b(v,l)+2\Delta))e^{(a(v,l)-2\Delta)r}e^{(b(v,l)-2\Delta)s}.
\end{multline*}
Combining the RHS and LHS with the density for $\tilde P^{(v,t)}(dr,ds)$, we
arrive at the inequality
\begin{multline*}
|\tilde P^{(v,t)}(dr,ds)-P^{(v,t)}(dr,ds)|  \\  \leq 
\left(1-e^{-2\Delta r}+\frac{2\Delta}{a(v,l)}\right)
 \left(e^{2\Delta s}-e^{-2\Delta s}+
 \frac{2\Delta (e^{2\Delta s}+e^{-2\Delta s})}{b(v,l)}\right)
  \tilde P^{(v,t)}(dr,ds).
\end{multline*}
Define $C(v,l,\Delta)=$
\begin{align}
\sup_{r,s<\Delta}&\left(1-e^{-2\Delta r}+\frac{2\Delta}{a(v,l)}\right)
 \left(e^{2\Delta s}-e^{-2\Delta s}+
 \frac{2\Delta e^{2\Delta s}+2\Delta e^{-2\Delta s}}{b(v,l)}\right) \\
&\leq \label{eq:defcvl}
 \left(1-e^{-2\Delta^2}+\frac{2\Delta}{a(-R,l)}\right)
\left(e^{2\Delta^2}-e^{-2\Delta^2}+
\frac{2\Delta e^{2\Delta^2}+2\Delta}{b(R,l)}\right).
\end{align}
It is clear that 
\[\displaystyle \lim_{\Delta,l\rightarrow 0}\sup_{|v|\leq R}C(v,l,\Delta)\rightarrow 0.\]

\end{rem}

\newcommand{\prti}{I}
\newcommand{\prtii}{I\!I}
\newcommand{\prtiii}{I\!I\!I}

\begin{lem}\label{lem:techfirst}
\[ ||a_l^{VV}(v,t)-2||\rightarrow 0. \]
\end{lem}
\begin{proof}
We first need to show that
\begin{align}\label{eq:tavv}
\frac{1}{l^2}\int_{E_l^c}(r-s)^2 \tilde P^{(v,t)}(dr,ds)\rightarrow 0.
\end{align}
Note that if $0\leq r\leq 1/2$ and $0 \leq s \leq 1/2$, then
$(r,s)\in E_l$ for $l<\sqrt{3}/2$.  Then
\begin{align*}
\frac{1}{l^2}\int_{E_l^c}& r^2 \tilde P^{(v,t)}(dr,ds) \leq
\frac{1}{l^2}\int\limits_{r>1/2}r^2 \tilde P^{(v,t)}(dr,ds)
  +\frac{1}{l^2}\int\limits_{s>1/2}r^2 \tilde P^{(v,t)}(dr,ds) \\
&= \frac{1}{l^2}\left(\int_{1/2}^{\infty}r^2 a(v,l)e^{-a(v,l)r}\ dr
   +E(\uls)^2 P(\dls>1/2)  \right) \\
&= \frac{1}{l^2}\left(\frac{2+a(v,l)+a^2(v,l)/4}{a^2(v,l)}e^{-a(v,l)/2}
   +\frac{2}{a^2(v,l)}e^{-b(v,l)/2}  \right) \\
&\leq \frac{1}{l^2}\left(\frac{2+a(R,l)+a^2(R,l)/4}{a^2(-R,l)}e^{-a(-R,l)/2)}
   +\frac{2}{a^2(-R,l)}e^{-b(R,l)/2}\right).
\end{align*}
Since $la(v,l)\rightarrow 1$ as $l\rightarrow 0$, we can
show that this last term will converge
to $0$ if we can show that $a^2(R,l)\exp(-a(-R,l)/2)\rightarrow 0$:
\begin{align*}
a^2(R,l)e^{-a(-R,l)/2} &\leq (1/l+2R)^2 e^{-1/l+2R},
\end{align*}
and we know that $l^{-m} e^{-1/l}\rightarrow 0$ as $l\rightarrow 0$ for all $m$.

The same procedure shows that
\begin{align}\label{eq:elc}
\frac{1}{l^2}\int_{E_l^c}s^2 \tilde P^{(v,t)}(dr,ds)\rightarrow 0.
\end{align}
and these calculations (noting that $(r-s)^2\leq r^2+s^2$) give \eqref{eq:tavv}.

Since $\uls$ and $\dls$ are independent exponential random variables,
\begin{align*}
\tilde{a}^{VV}_l(v,t) &=
\frac{1}{l^2}\int_{E_l}(r-s)^2 \tilde P^{(v,t)}(dr,ds) \\
&=\frac{1}{l^2}E (\uls-\dls)^2 
    -\frac{1}{l^2}\int_{E_l^c}(r-s)^2 \tilde P^{(v,t)}(dr,ds)\\
&=\frac{1}{l^2}
  \left(\frac{2}{a^2(v,l)}-\frac{2}{a(v,l)b(v,l)}+\frac{2}{b^2(v,l)}\right)
   -\prti.
\end{align*}
Since $a(v,l)$ and $b(v,l)$ are monotone, and $l a(v,l)$ and $l
b(v,l)$ converge to $1$ for all $v$, and $\prti\rightarrow 0,$ the RHS of
the last equation converges to $2$ uniformly for $|v|\leq R$.  Combining
this with equation \eqref{eq:tavv} shows that
\[ ||\tilde a_l^{VV}-2||\rightarrow 0. \]

Next we compare $a_l^{VV}$ and $\tilde a_l^{VV}$.
\begin{align*}
a_l^{VV}&-\tilde{a}_l^{VV} = 
 \frac{1}{l^2}\int_{E_l}(r-s)^2 
  (P^{(v,t)}(\ul\in dr,\dl\in ds)-P^{(v,t)}(\uls\in dr,\dls\in ds)) \\
&= \frac{1}{l^2}\int\limits_{E_l\cap\{r+s<\sqrt{l}\}}(r-s)^2 
  (P^{(v,t)}(\ul\in dr,\dl\in ds)-P^{(v,t)}(\uls\in dr,\dls\in ds))\\
  &+\frac{1}{l^2}\int\limits_{E_l\cap\{r+s>\sqrt{l}\}}(r-s)^2 
  (P^{(v,t)}(\ul\in dr,\dl\in ds)-P^{(v,t)}(\uls\in dr,\dls\in ds)) \\
&= \prtii+\prtiii.
\end{align*}
Since $(r-s)^2\leq 1$ on $E_l$, the second term $\prtiii$ can be bounded as follows:
\begin{gather}
\begin{split}\label{eq:dcalc}
\prtiii &\leq \frac{1}{l^2}[P(\ul>\sqrt{l}/2)+P(\dl>\sqrt{l}/2,\ul<\sqrt{l}/2)] \\
 &\qquad+ \frac{1}{l^2}[P(\uls>\sqrt{l}/2)+P(\dls>\sqrt{l}/2)] 
\end{split} \\
\leq \frac{1}{l^2}[e^{-\sqrt{l}a(v,l)/2}+e^{-\sqrt{l}(b(v,l)-2K\sqrt{l})/2}+
  e^{-\sqrt{l}a(v,l)/2}+e^{-\sqrt{l}b(v,l)/2}] \label{Dineq}
\end{gather}
since if $\ul<\sqrt{l},$ the rate of $\dl<b(v,l)+2K\sqrt{l}$.
We can compute that
\begin{align*}
\frac{1}{l^2}e^{-\sqrt{l}a(v,l)/2} &<
 \frac{1}{l^2}e^{-\sqrt{l}(1/l-2R)} \\
 &=\frac{1}{l^2}e^{-1/\sqrt{l}}e^{2R\sqrt{l}}.
\end{align*}
The term $\exp(-1/\sqrt l)/l^2\rightarrow 0$ and 
$\exp(2R\sqrt l)\rightarrow 1$.

To bound the term $\prtii$, we observe that on the set
$E_l\cap\{r+s<\sqrt{l}\},$ $\vtw$ must lie in the interval 
$(v-K\sqrt{l},v+K\sqrt{l})$.  Then
\begin{align*}
 \prtii&\leq C(v,l,\sqrt{l})
   \frac{1}{l^2}\int_{{E_l\cap\{r+s<\sqrt{l}\}}}(r-s)^2\tilde P^{(v,t)}(dr,ds) \\
 &\leq C(v,l,\sqrt{l})
   \frac{1}{l^2}\int_{r,s>0}(r-s)^2\tilde P^{(v,t)}(dr,ds) \\
 &=C(v,l,\sqrt{l})
  \frac{1}{l^2}
  \left(\frac{2}{a^2(v,l)}-\frac{2}{a(v,l)b(v,l)}+\frac{2}{b^2(v,l)}\right).
\end{align*}
As before, $(2/a^2-2/ab+2/b)/l^2\rightarrow 2$,
and as shown above, $C(v,l,\sqrt{l})\rightarrow 0$ uniformly on
compact sets as $l\rightarrow 0$, so we conclude that
$||a_l^{VV}(v,t)-\tilde{a}_l^{VV}(v,t)||\rightarrow 0$.  Finally, we
use the triangle inequality to conclude that
$||a_l^{VV}(v,t)-2||\rightarrow 0$.
\end{proof}

\newcommand{\dcr}{{(\ref{eq:dcalc}) on page \pageref{eq:dcalc}}}


\begin{lem}\label{lem:btconv}
\[ ||b_l^T(v,t)-2||\rightarrow 0. \]
\end{lem}
\begin{proof}
We will first show that
\begin{align}\label{eq:btl}
\frac{1}{l}\int_{E_l^c}(r+s) \tilde P^{(v,t)}(dr,ds)\rightarrow 0.
\end{align}
As in the previous case,
\begin{align*}
\frac{1}{l}\int_{E_l^c}& r \tilde P^{(v,t)}(dr,ds) \leq
\frac{1}{l}\int\limits_{r>1/2}r \tilde P^{(v,t)}(dr,ds)
  +\frac{1}{l}\int\limits_{s>1/2}r \tilde P^{(v,t)}(dr,ds) \\
&= \frac{1}{l}\left(\int_{1/2}^{\infty}r a(v,l)e^{-a(v,l)r}\ dr
   +E(\uls) P(\dls>1/2)  \right) \\
&= \frac{1}{l}\left(\frac{1+a(v,l)/2}{a(v,l)}e^{-a(v,l)/2}
   +\frac{1}{a(v,l)}e^{-b(v,l)/2}  \right) \\
&\leq \frac{1}{l}\left(\frac{1+a(R,l)/2}{a(-R,l)}e^{-a(-R,l)/2}
   +\frac{1}{a(R,l)}e^{-b(R,l)/2}  \right). \\
\end{align*}
Since $la(v,l)\rightarrow 1$ as $l\rightarrow 0$, we can show that this last term will converge
to $0$ by showing that $a(R,l)\exp(-a(-R,l)/2)\rightarrow 0$:
\begin{align*}
a(R,l)e^{-a(-R,l)/2} &\leq (1/l+2R) e^{-1/l+2R},
\end{align*}
and we know that $l^{-m} e^{-1/l}\rightarrow 0$ as $l\rightarrow 0$ for all $m$.

The same procedure shows that
\begin{align*}
\frac{1}{l}\int_{E_l^c}s \tilde P^{(v,t)}(dr,ds)\rightarrow 0.
\end{align*}
and combining them gives \eqref{eq:btl}.

Since $\uls$ and $\dls$ are independent exponential random variables,
\begin{align*}
\tilde{b}^T_l(v,t) &=
\frac{1}{l}\int_{E_l}(r+s) \tilde P^{(v,t)}(dr,ds) \\
&=\frac{1}{l}E (\uls+\dls) 
    -\frac{1}{l}\int_{E_l^c}(r+s) \tilde P^{(v,t)}(dr,ds)\\
&=\frac{1}{l}
  \left(\frac{1}{a(v,l)}+\frac{1}{b(v,l)}\right)
   -\prti.
\end{align*}
Since $a(v,l)$ and $b(v,l)$ are monotone, and $l a(v,l)$ and $l b(v,l)$ converge
to $1$ for all $v$, and $\prti\rightarrow 0,$ this quantity converges 
to $2$ uniformly for $|v|\leq R$.

Next we compare $b_l^T$ and $\tilde b_l^T$.
\begin{align*}
b_l^{T}&-\tilde{b}_l^{T} = 
 \frac{1}{l}\int_{E_l}(r+s) 
  (P^{(v,t)}(\ul\in dr,\dl\in ds)-P^{(v,t)}(\uls\in dr,\dls\in ds)) \\
&= \frac{1}{l}\int\limits_{E_l\cap\{r+s<\sqrt{l}\}}(r+s) 
  (P^{(v,t)}(\ul\in dr,\dl\in ds)-P^{(v,t)}(\uls\in dr,\dls\in ds))
\\&+
  \frac{1}{l}\int\limits_{E_l\cap\{r+s>\sqrt{l}\}}(r+s) 
  (P^{(v,t)}(\ul\in dr,\dl\in ds)-P^{(v,t)}(\uls\in dr,\dls\in ds)) \\
&= \prtii+\prtiii.
\end{align*}
Since $r+s\leq 1/l$ on $E_l$, the second term $\prtiii$ can be bounded as follows:
\begin{align*}
\prtiii &\leq \frac{1}{l^2}
  [P(\ul>\sqrt{l}/2)+P(\dl>\sqrt{l}/2,\ul<\sqrt{l}/2)] \\
 &\qquad+ \frac{1}{l^2}[P(\uls>\sqrt{l}/2)+P(\dls>\sqrt{l}/2)]
\end{align*}
which we know converges to $0$, since it is exactly the same as \dcr.

To bound the term $\prtii$, we observe that on the set
$E_l\cap\{r+s<\sqrt{l}\},$ $\vtw$ must lie in the interval 
$(v-K\sqrt{l},v+K\sqrt{l})$.  Then
\begin{align*}
 \prtii&\leq C(v,l,\sqrt{l})
   \frac{1}{l}\int_{{E_l\cap\{r+s<\sqrt{l}\}}}(r+s)\tilde P^{(v,t)}(dr,ds) \\
 &\leq C(v,l,\sqrt{l})
   \frac{1}{l}\int_{r,s>0}(r+s)\tilde P^{(v,t)}(dr,ds) \\
 &=C(v,l,\sqrt{l})
  \frac{1}{l}
  \left(\frac{1}{a(v,l)}+\frac{1}{b(v,l)}\right).
\end{align*}
As before, $(1/a+1/b)/l\rightarrow 2$,
and as shown above, $C(v,l,\sqrt{l})\rightarrow 0$ uniformly on
compact sets as $l\rightarrow 0$, so we conclude that
$||b_l^{T}(v,t)-\tilde{b}_l^{T}(v,t)||\rightarrow 0$.  Finally, we use
the triangle inequality to conclude that
$||b_l^{T}(v,t)-2||\rightarrow 0$.
\end{proof}


\begin{lem}
\[ ||b_l^V(v,t)+2v||\rightarrow 0 \]
\end{lem}
\begin{proof}
We first show that
\begin{align}\label{eq:bvl}
\frac{1}{l^2}\int_{E_l^c}(r-s) \tilde P^{(v,t)}(dr,ds)\rightarrow 0.
\end{align}
As in the previous case,
\begin{align*}
\frac{1}{l^2}\int_{E_l^c}& r \tilde P^{(v,t)}(dr,ds) \leq
\frac{1}{l^2}\int\limits_{r>1/2}r \tilde P^{(v,t)}(dr,ds)
  +\frac{1}{l^2}\int\limits_{s>1/2}r \tilde P^{(v,t)}(dr,ds) \\
&\leq \frac{1}{l^2}\left(\frac{1+a(R,l)/2}{a(-R,l)}e^{-a(-R,l)/2}
   +\frac{1}{a(R,l)}e^{-b(R,l)/2}  \right).
\end{align*}
Because $la(v,l)\rightarrow 1$ as $l\rightarrow 0$, 
this last term will converge
to $0$ if we can show that $(1/l)a(R,l)\exp(-a(-R,l)/2)\rightarrow 0$:
\begin{align*}
(1/l)a(R,l)e^{-a(-R,l)/2} &\leq (1/l^2+2R/l) e^{-1/l+2R},
\end{align*}
and we know that $l^{-m} e^{-1/l}\rightarrow 0$ as $l\rightarrow 0$ for all $m$.

The same procedure shows that
\begin{align*}
\frac{1}{l^2}\int_{E_l^c}s \tilde P^{(v,t)}(dr,ds)\rightarrow 0.
\end{align*}
and combining them gives \eqref{eq:bvl}.

Since $\uls$ and $\dls$ are independent exponential random variables,
\begin{align*}
\tilde{b}^V_l(v,t)&=
\frac{1}{l^2}\int_{E_l}(r-s) \tilde P^{(v,t)}(dr,ds) \\
&=\frac{1}{l^2}E (\uls-\dls) 
    -\frac{1}{l^2}\int_{E_l^c}(r-s) \tilde P^{(v,t)}(dr,ds)\\
&=\frac{1}{l^2}
  \left(\frac{1}{a(v,l)}-\frac{1}{b(v,l)}\right)
   -\prti \\
&=\frac{1}{l^2}
  \frac{1}{a(v,l)b(v,l)}(b(v,l)-a(v,l))
   -\prti \\
&=\frac{1}{l^2}
  \frac{1}{a(v,l)b(v,l)}\frac{v(e^{-vl}-e^{vl})}{\sinh(vl)}
   -\prti \\
&=\frac{1}{l^2}
  \frac{1}{a(v,l)b(v,l)}(-2v)
   -\prti.
\end{align*}
Since $a(v,l)$ and $b(v,l)$ are monotone, and $l a(v,l)$ and $l b(v,l)$ converge
to $1$ for all $v$, and $\prti\rightarrow 0,$ the last equation converges 
to $-2v$ uniformly for $|v|\leq R$.

Next we compare $b_l^V$ and $\tilde b_l^V$.
\begin{align*}
|b_l^{V}&-\tilde{b}_l^{V}| \leq 
 \frac{1}{l^2}\int_{E_l}|r-s| 
  (P^{(v,t)}(\ul\in dr,\dl\in ds)-P^{(v,t)}(\uls\in dr,\dls\in ds)) \\
&= \frac{1}{l^2}\int_{E_l\cap\{r+s<\sqrt{l}\}}|r-s| 
  (P^{(v,t)}(\ul\in dr,\dl\in ds)-P^{(v,t)}(\uls\in dr,\dls\in ds)) \\
&+
  \frac{1}{l^2}\int\limits_{E_l\cap\{r+s>\sqrt{l}\}}
  |r-s| 
  (P^{(v,t)}(\ul\in dr,\dl\in ds)-P^{(v,t)}(\uls\in dr,\dls\in ds)) \\
&= \prtii+\prtiii.
\end{align*}
Since $|r-s|\leq 1$ on $E_l$, the second term $\prtiii$ can be bounded as follows:
\begin{align*}
\prtiii &\leq \frac{1}{l^2}
  [P(\ul>\sqrt{l}/2)+P(\dl>\sqrt{l}/2,\ul<\sqrt{l}/2)] \\
 &\qquad+ \frac{1}{l^2}[P(\uls>\sqrt{l}/2)+P(\dls>\sqrt{l}/2)]
\end{align*}
which converges to $0$, as \dcr.

To bound the term $\prtii$, we observe that on the set
$E_l\cap\{r+s<\sqrt{l}\},$ $\vtw$ must lie in the interval 
$(v-K\sqrt{l},v+K\sqrt{l})$.  Then
\begin{align*}
 \prtii&\leq \frac{1}{l}C(v,l,\sqrt{l})
   \frac{1}{l}\int_{{E_l\cap\{r+s<\sqrt{l}\}}}(r+s)\tilde P^{(v,t)}(dr,ds) \\
 &\leq \frac{1}{l}C(v,l,\sqrt{l})
   \frac{1}{l}\int_{r,s>0}(r+s)\tilde P^{(v,t)}(dr,ds) \\
 &=\frac{1}{l}C(v,l,\sqrt{l})
  \frac{1}{l}
  \left(\frac{1}{a(v,l)}+\frac{1}{b(v,l)}\right).
\end{align*}
As before, $(1/a+1/b)/l\rightarrow 2$.
For this case, we need to show that 
$C(v,l,\sqrt{l})/l\rightarrow 0$
uniformly on compact sets as
$l\rightarrow 0$.
From \eqref{eq:defcvl}, it is enough to show that
\begin{equation*}
\frac 1 l\left(1-e^{-2l}+\frac{2\sqrt{l}}{a(-R,l)}\right)\rightarrow 2
\end{equation*}
(in fact, any finite limit will do), and this is an application of
l'H\^opital's rule.  We conclude that
$||b_l^{V}(v,t)-\tilde{b}_l^{V}(v,t)||\rightarrow 0$.  Finally, we use
the triangle inequality to conclude that
$||b_l^{V}(v,t)+2v||\rightarrow 0$.
\end{proof}


\raggedbottom
\begin{lem}\label{lem:attzero}
\[ ||a_l^{TT}(v,t)||\rightarrow 0. \]
\end{lem}

%
\begin{proof}
Since $\uls$ and $\dls$ are independent exponential random variables,
\begin{align*}
\tilde{a}_l^{TT}(v,t)&=
\int_{E_l}(r+s)^2 \tilde P^{(v,t)}(dr,ds) \\
&\leq E(\uls+\dls)^2 \\
&=\frac{2}{a^2(v,l)} +\frac{2}{a(v,l)b(v,l)}
    +\frac{2}{b^2(v,l)}.
\end{align*}
Since $a(v,l)$ and $b(v,l)$ are monotone, and $a(v,l)$ and $b(v,l)$ converge
to $\infty$ for all $v$, the last equation converges 
to $0$ uniformly for $|v|\leq R$.

Next we compare $a_l^{TT}$ and $\tilde a_l^{TT}$.
\begin{align*}
|a_l^{TT}&-\tilde{a}_l^{TT}| \leq 
 \int_{E_l}(r+s)^2 
  (P^{(v,t)}(\ul\in dr,\dl\in ds)-P^{(v,t)}(\uls\in dr,\dls\in ds)) \\
&= \int_{E_l\cap\{r+s<\sqrt{l}\}}(r+s)^2 
  (P^{(v,t)}(\ul\in dr,\dl\in ds)-P^{(v,t)}(\uls\in dr,\dls\in ds)) \\ 
&+
  \int\limits_{E_l\cap\{r+s>\sqrt{l}\}}(r+s)^2 
  (P^{(v,t)}(\ul\in dr,\dl\in ds)-P^{(v,t)}(\uls\in dr,\dls\in ds)) \\
&= \prti+\prtii.
\end{align*}
Since $(r+s)^2\leq 1/l^2$ on $E_l$, the second term $\prtii$ can be bounded as follows:
\begin{align*}
\prtii &\leq \frac{1}{l^2}
  [P(\ul>\sqrt{l}/2)+P(\dl>\sqrt{l}/2,\ul<\sqrt{l}/2)] \\
 &\qquad+ \frac{1}{l^2}[P(\uls>\sqrt{l}/2)+P(\dls>\sqrt{l}/2)]
\end{align*}
which converges to $0$, since it is the same as \dcr.

To bound the term $\prti$, we observe that on the set
$E_l\cap\{r+s<\sqrt{l}\},$ $\vtw$ must lie in the interval 
$(v-K\sqrt{l},v+K\sqrt{l})$.  Then
\begin{align*}
 \prti&\leq C(v,l,\sqrt{l})
   \int_{{E_l\cap\{r+s<\sqrt{l}\}}}(r+s)^2\tilde P^{(v,t)}(dr,ds) \\
 &\leq C(v,l,\sqrt{l})
   \int_{r,s>0}(r+s)^2\tilde P^{(v,t)}(dr,ds) \\
 &=C(v,l,\sqrt{l})
  \left(\frac{2}{a^2(v,l)}+\frac{2}{a(v,l)b(v,l)}+\frac{2}{b^2(v,l)}\right).
\end{align*}
In this case, both terms converge to $0$ as $l\rightarrow 0$.
We conclude that 
$||a_l^{TT}(v,t)-\tilde{a}_l^{TT}(v,t)||\rightarrow 0$.  
Finally, we conclude that $||a_l^{TT}(v,t)||\rightarrow 0$.
\end{proof}


\begin{lem}
\[ ||a_l^{VT}(v,t)||\rightarrow 0. \]
\end{lem}
%
\begin{proof}
Since $\uls$ and $\dls$ are independent exponential random variables,
\begin{align*}
\tilde{a}^{VT}_l(v,t) &=
\frac 1 l \int_{E_l}(r+s)(r-s) \tilde P^{(v,t)}(dr,ds) \\
&\leq \frac 1 l E[(\uls)^2+(\dls)^2] \\
&=\frac 1 l \left(\frac{2}{a^2(v,l)}+\frac{2}{b^2(v,l)}\right).
\end{align*}
Since $a(v,l)$ and $b(v,l)$ are monotone, $la(v,l)$ and $lb(b,l)$ converge
to $1$,
and $a(v,l)$ and $b(v,l)$ converge
to $\infty$ for all $v$, the last equation converges 
to $0$ uniformly for $|v|\leq R$.

Next we compare $a_l^{VT}$ and $\tilde a_l^{VT}$.
\begin{align*}
|a_l^{VT}&-\tilde{a}_l^{VT}| \leq 
 \frac 1 l \int\limits_{E_l}|r^2-s^2| 
  (P^{(v,t)}(\ul\in dr,\dl\in ds)-P^{(v,t)}(\uls\in dr,\dls\in ds)) \\
&= \frac 1 l \int\limits_{E_l\cap\{r+s<\sqrt{l}\}}|r^2-s^2| 
  (P^{(v,t)}(\ul\in dr,\dl\in ds)-P^{(v,t)}(\uls\in dr,\dls\in ds))\\
&+
  \frac 1 l \int\limits_{E_l\cap\{r+s>\sqrt{l}\}}|r^2-s^2| 
  (P^{(v,t)}(\ul\in dr,\dl\in ds)-P^{(v,t)}(\uls\in dr,\dls\in ds)) \\
&= \prti+\prtii.
\end{align*}
Since $|r^2-s^2|\leq 1/l$ on $E_l$, the second term $\prtii$ can be bounded as follows:
\begin{align*}
\prtii &\leq \frac{1}{l^2}
  [P(\ul>\sqrt{l}/2)+P(\dl>\sqrt{l}/2,\ul<\sqrt{l}/2)] \\
 &\qquad+ \frac{1}{l^2}[P(\uls>\sqrt{l}/2)+P(\dls>\sqrt{l}/2)]
\end{align*}
which we know converges to $0$, since it is exactly the same as \dcr.

To bound the term $\prti$, we observe that on the set
$E_l\cap\{r+s<\sqrt{l}\},$ $\vtw$ must lie in the interval 
$(v-K\sqrt{l},v+K\sqrt{l})$.  Then
\begin{align*}
 \prti&\leq C(v,l,\sqrt{l})
   \int_{{E_l\cap\{r+s<\sqrt{l}\}}}(r+s)^2\tilde P^{(v,t)}(dr,ds) \\
 &\leq C(v,l,\sqrt{l})
   \int_{r,s>0}(r+s)^2\tilde P^{(v,t)}(dr,ds) \\
 &=C(v,l,\sqrt{l})
  \left(\frac{2}{a^2(v,l)}+\frac{2}{a(v,l)b(v,l)}+\frac{2}{b^2(v,l)}\right).
\end{align*}
In this case, both terms converge to $0$ as $l\rightarrow 0$, so
$||a_l^{VT}(v,t)-\tilde{a}_l^{VT}(v,t)||\rightarrow 0$.  
Finally, we conclude by the triangle inequality that 
$||a_l^{VT}(v,t)||\rightarrow 0$.
\end{proof}

\begin{lem}\label{lem:techlast}
\[ \frac{1}{l^2}
 \Delta_l^{\varepsilon}(v,t)\rightarrow 0 \]
uniformly for $|v|<R$.  
\end{lem}
\begin{proof}
In fact, we need only observe that for $l<1/2$,
the set $\{(u,v) : u+v<\sqrt l\} \subset E_l$.  Then
\begin{align*}
\Delta_l^{\varepsilon}(v,t) &\leq \frac{1}{l^2}
  [P(\ul>\sqrt{l}/2)+P(\dl>\sqrt{l}/2,\ul<\sqrt{l}/2)],
\end{align*}
which we have shown to converge above (again, \dcr).
\end{proof}

\begin{proof}[Proof of Lemma \ref{lem:genconv}]
We need to show that the martingale problem for
\[ L= \frac{\partial^2}{\partial v^2}-2v\frac{\partial}{\partial v}
+2\frac{\partial}{\partial T} \]
has a unique solution.  Since the coefficients are either bounded or linear, we
can use Theorem 5.2.9 of \cite{ks:91}.  Once uniqueness of the solution of the
martingale problem for $L$ is established, Lemma \ref{lem:genconv} follows directly 
from Theorem 11.2.3 of \cite{sv} and Lemmas \ref{lem:techfirst}--\ref{lem:techlast}.
\end{proof}


An issue with Lemma \ref{lem:genconv} is that the measures $P^l_{v,t}$
are associated with $\hat V^l(t)$ rather than with the original
$\tilde V^l(\tau)$.  In fact, the convergence holds for $\tilde V^l(\tau)$ 
as well.

\begin{thm}
The process $\tilde V^l(\tau)$ converges weakly to the Ornstein-Uhlenbeck 
process.
\end{thm}
\begin{proof}
  By the symmetry of $a$ and $b$, the process interpolated on the
  other side (along $\sigma_j$'s) has the same limit as $\hat V^l(t)$.
Since the processes interpolated along
  the top and the bottom of the sawtooth process converge to the same
  process, the whole sawtooth process must converge if we show that the
distance between them is $0$, or equivalently, that the distance between
the two processes converges to $0$ uniformly on finite time intervals.
This follows from Lemmas \ref{lem:btconv} and \ref{lem:attzero}.

The construction of $\hat V^l(t)$ involved a time change because we forced
each switch of the process to have duration $l^2$.  The term 
$2\frac{\partial}{\partial T}$ of the generator $L$ indicates that in the
limit, this is twice the duration of the actual time between switches.
We can restore the original clock by dividing the generator by two.  When
we do so, we find that the spatial component of the process has generator
corresponding to solutions of the SDE
\[ dX_t=dB_t-X_t\ dt. \]
This is an Ornstein-Uhlenbeck process.
\end{proof}


\section{A Pair of Brownian Motions Separated by an Inert Particle}\label{ch:twobm}

\newcommand{\bxk}{\mathcal{X}^K}
\newcommand{\byk}{\mathcal{Y}^K}
\newcommand{\xki}{X_1}
\newcommand{\xkii}{X_2}
\newcommand{\yk}{Y}
\newcommand{\lki}{L_1}
\newcommand{\lkii}{L_2}
\newcommand{\vk}{V}
\newcommand{\tk}{T}
\newcommand{\defn}{\equiv}


\newcommand{\tinf}{T_\infty}

In this section, we consider an arrangement of two Brownian particles
$\xki$ and $\xkii$ separated by an inert particle $\yk$ in $\Re$.  
More precisely, we
construct processes $\xki(t)\leq\yk(t)\leq\xkii(t)$,
where the interactions between $\xki(t)$ and $\yk(t)$ and between 
$\yk(t)$ and $\xkii(t)$ are as described in Section \ref{ch:sk}.

\begin{figure}[ht]
\centering
\includegraphics[width=4in]{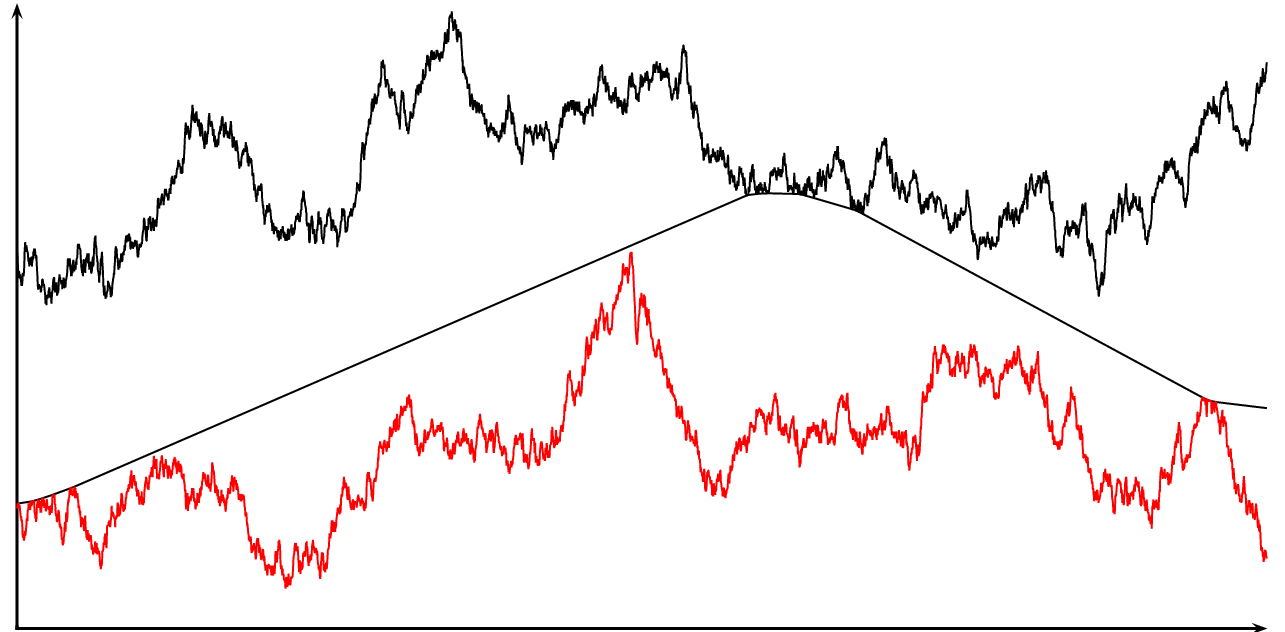}
\caption{Two Brownian particles separated by an inert particle}
\end{figure}

A method of construction different from that in Section \ref{ch:sk}
is needed if the two Brownian particles ever meet.  Instead
we introduce a random variable $\tinf$ to represent the first
meeting of the two Brownian particles $\xki(t)$ and $\xkii(t)$.  
In fact, we will show that with probability one, $\tinf=\infty$.

\begin{thm}\label{thm:twopart}
Given independent Brownian motions $B_1$ and $B_2$, 
with $B_j(0)=0$,
constants $x>0$, $0\leq y\leq x$, $v\in\Re$, and $K>0$, there exist 
unique processes
$\lki(t)$ and $\lkii(t)$, and a random time $\tinf$, 
satisfying the following conditions:
\begin{enumerate}
\item $\xki(t)\leq\yk(t)\leq\xkii(t),\ 0\leq t\leq \tinf,$ where
\begin{enumerate}
\item $\xki(t)\defn B_1(t)-\lki(t)$,
\item $\xkii(t)\defn x+B_2(t)+\lkii(t)$,
\item $\displaystyle\vk(t)\defn v+K(\lki(t)-\lkii(t))$, and
\item $\displaystyle\yk(t)\defn y+\int_0^t\vk(s)\,ds$,
\end{enumerate}
\item $\lki(t)$ and $\lkii(t)$ are continuous, nondecreasing functions
with $\lki(0)=\lkii(0)=0$,
\item $\lki(t)$ and $\lkii(t)$ are flat off the sets
  $\{t:\xki(t)=\yk(t)\}$ and $\{t: \xkii(t)=\yk(t)\}$, resp.
\item $\tinf=\inf\{t: \xki(t)=\xkii(t)\}$.
\end{enumerate}
\end{thm}
\begin{proof}
The construction method in the first section and a sequence of stopping
times can be used to construct this process up to the stopping time $\tinf$, the limit of the stopping times used in the construction.
After time $\tinf$ the process is not well-defined, but we show below that 
$P(\tinf=\infty)=1$.
\end{proof}

We define $\bxk_x(t)\defn (\xki(t),\xkii(t),\yk(t), \vk(t))$ for the
processes constructed with initial state $y=0$, $v=0$, and constant
$K$.  The following lemma describes the scaling law of the process.

\begin{lem}\label{lem:scale}
$\varepsilon\bxk_x(t/\varepsilon^2)\overset{d}{=}\mathcal{X}^{K/\varepsilon}_{\varepsilon
    x}(t)$.
\end{lem}
\begin{proof}
By Brownian scaling, the $\xki$ and $\xkii$ components remain Brownian motions,
and by uniqueness of local time, $\lki$ and $\lkii$ have the same scaling.
However, by the chain rule, the slope of the $\yk$ component has been multiplied
by $1/\varepsilon$ for each $t$.
\end{proof}

The rest of the section concerns the proof that $\tinf=\infty$ \as
\newcommand{\xrt}{X^{\rho,T}}
\newcommand{\lrt}{L^{\rho,T}}
\newcommand{\lrti}{\lrt_\infty}
\begin{thm}\label{thm:vnotbig}
Define a process $\xrt(t)$ for $T>0$ and $\rho\in(0,1)$ as follows.  Let
\[ \Delta(t)=
	\begin{cases} 
		0 & t<T \\ 
		-\rho L(T)(t-T) & t\geq T
	\end{cases} 
\]
By previous results, there are  unique $\lrt(t)$ and $L(t)$ such that 
\begin{align*} 
\xrt(t)&=[B(t)+\Delta(t)]+\lrt(t)+\int_0^t \lrt(s)\ ds\geq 0, \\
 X(t)&=B(t)+L(t)+\int_0^t L(s)\ ds\geq 0,
 \end{align*}
where $\lrt(t)$ and $L(t)$ are the local times of $\xrt(t)$ and $X(t)$, respectively, at zero.
Define $\lrti=\lim\limits_{t\rightarrow\infty}\lrt(t)$ and $L_\infty=\lim\limits_{t\rightarrow\infty}L(t)$.
Then 
\[ P(\lrti>l) \leq P(L_\infty>l)=\exp(-l^2). \]
\end{thm}
\begin{proof}
First note that $X(t)=\xrt(t)$ and $L(t)=\lrt(t)$ for $t\leq T$.  Also note that the drift of term of $\xrt(t)$ at
time $T$ is $(1-\rho)L(T)>0$.  After time $T$, $\xrt(t)$ may or may not return to the origin.  If not, then $X(t)$ also would not have returned to the origin ($B(t)\geq B(t)+\Delta(t)$), so $\lrti=(1-\rho)L_\infty$.

\newcommand{\trt}{\tau^{\rho,T}}
\newcommand{\bp}{\tilde{B}}
\newcommand{\lp}{\tilde{L}}
\newcommand{\xp}{\tilde{X}}

Otherwise, $\xrt$ returns to the origin at some time $\trt$.  Define 
 \[ S=\inf\{t \mid L(t)=(1-\rho)L(T) \} \]
Notice that $X(S)=0$ with probability $1$.  Construct a Brownian motion $\bp$ by 
deleting the time interval $(S,\trt)$ from $B(t)$:
\[ \bp(t)=\begin{cases}B(t) & t\leq S \\ B(t-S+\trt)-B(\trt)+B(S) &  t \geq S \end{cases}, \]
and an associated local time:
\[ \lp(t)=\begin{cases}L(t) & t\leq S \\ \lrt(t-S+\trt)-\rho L(T) &  t \geq S \end{cases} \]
and the associated reflected process with drift:
\[ \xp(t)=\begin{cases}X(t) & t\leq S \\ \xrt(t-S+\trt) &  t \geq S \end{cases}. \]
Note that $\bp(t)$ is a Brownian motion because $\trt$ is a stopping time and $S$ is depends
only on $B[S,T]$ and so is independent of $B[0,S]$.

We will show that $\xp(t)=\bp(t)+\lp(t)+\int_0^t\lp(s)\ ds$.  In fact, because of the pathwise uniqueness of solutions $L(t)$, we only need to check that $\bp,$ $\xp$, and $\lp$ are continuous at $S$.
\begin{align*}
\lp(S-) &= L(S)=(1-\rho)L(S) \\
\lp(S+) &= \lrt(\trt)-\rho L(T) \\
 &= \lrt(T)-\rho L(T) \\
 &= L(T)-\rho L(T) \\
 &= L(S-)
\end{align*}


The limit of $\lp(t)$ as $t\rightarrow\infty$ is $\lrti$ (pathwise).  But the limit of $\lp(t)$ will have the same distribution as $L_\infty$ because $\bp(t)$ is a Brownian motion. Since we have either decreased $L_\infty$ by a factor of $\rho$ or replaced it with a new copy with identical distribution, the inequality holds.
\end{proof}

\begin{thm}\label{thm:tinfty}
For the process constructed in Theorem \ref{thm:twopart}, 
$P(T=\infty)=1$.
\end{thm}
\begin{proof}
By the previous lemma, we may assume that $K=1$. We also assume that $v=0$.
We use slightly simplified versions of the process $X_1$ and $X_2$ below, which
incorporate the drift term ($Y(t)$ in the definiton), and which otherwise agree until time $T$ with
the definitions in Theorem \ref{thm:twopart}.

\[ X_1(t)=B_1(t)+\lki(t)+\int_0^t V(s)\ ds\geq 0 \]
\[ X_2(t)=B_2(t)-\lkii(t)+\int_0^t V(s)\ ds\leq 0 , \]
where $\lki(t),$ $\lkii(t)$ are the local times of $X_1(t)$ and $X_2(t)$ at the origin, $B_1(0)>0$ and $B_2(0)=0,$ and
\[ V(t)=\begin{cases} \lki(t)-\lkii(t) & t<T \\ 0 & t\geq T\end{cases}, \]
with $T$ a stopping time defined below.

Define $T_0=0$, $T_{j+1}=\inf\{t>T_j \mid V(t)=0 \}$, and define $\tinf=\lim T_j$.
On any of the intervals $[T_j,T_{j+1}]$ (say that $X_1(T_j)=0$), the term
$V(t)$ behaves exactly as in the case of one Brownian particle and one inert particle, except that $V(t)$ may decrease when $X_2(t)=0$.  So $V(t)$ is dominated in distribution by $L_\infty$.

Using the previous theorem, we can check Novikov's condition and then apply the Girsanov theorem to $(X_1, X_2)$:
\begin{align*}
 E\exp\left(\frac{1}{2}\int_j^{j+1} (V(t))^2 dt  \right)
  &\leq E\exp\left(\frac{1}{2}\int_j^{j+1} (L_\infty)^2 dt  \right)\\
  &\leq E\exp\left(\frac{1}{2}(L_\infty)^2 dt  \right) \\
  &= \int_0^\infty \exp\left(\frac{1}{2}s^2\right)P(L_\infty\in ds)\\
  &= \int_0^\infty\left[\int_0^s \exp\left(\frac{1}{2}r^2\right)r\ dr+1\right]P(L_\infty\in ds) \\
  &= \int_0^\infty\exp\left(\frac{1}{2}r^2\right)r\int_r^\infty P(L_\infty\in ds)\ dr+1 \\
  &= \int_0^\infty\exp\left(\frac{1}{2}r^2\right)r\exp(-r^2)dr+1\\
  &<\infty.
\end{align*}
We can now apply Girsanov to see that under some measure, $(X_1,X_2)$ is a standard reflected Brownian motion in the quadrant $\{(x,y)\mid x>0, y>0\}$.  Observe that
if $X_1(T_j)=0$, then $X_2(T_{j+1})=0$ and $X_1(T_{j+2}=0$.  Then $T_\infty<\infty$ implies that the reflected Brownian motion hits the origin, an event with probability zero.  Therefore, $P(T_\infty=\infty)=1$.
\end{proof}

\subsection{The limiting process is Bess(2)}

\newcommand{\bx}{{|(y-x,t)|<1}}
\newcommand{\pik}{\Pi_K(x,dy,dt)}

In this section, we wish to determine the law of the process described
in Theorem \ref{thm:twopart} as the constant $K\rightarrow\infty$.
As in the previous section, we approach the limit distribution through a 
Markov chain.  We introduce the stopping times $T_j$ defined by
\[ T_0=0,\ T_{j+1}=\inf\{t>T_j: V(t)=0\}, \]
and a two-dimensional Markov chain $\{\byk(j)\}_{j=0}^\infty$
defined by 
\[ \byk(j)=(\xkii(T_j)-\xki(T_j),T_j).\]
We denote the transition probabilities of $\byk(j)$ by
\[ \pik=P\left(\byk(j+1)\in(dy,dt)\mid\byk(j)=(x,0)\right), \]
noting that $T_{j+1}-T_j$ is independent of the value of $T_j$.





Now that our processes are defined, we focus on the transition probabilities
of $\{\byk_j\}_j$.  The following definitions correspond to those in
\cite[section 11.2]{sv}, with $h=1/\sqrt{K}$.
\begin{align*}
 b_K^X(x) &= \sqrt{K}\int_\bx (y-x)\ \pik, \\
 b_K^T(x) &= \sqrt{K}\int_\bx t\ \pik, \\
 a_K^{XX}(x) &= \sqrt{K}\int_\bx (y-x)^2\ \pik, \\
 a_K^{XT}(x) &= \sqrt{K}\int_\bx (y-x)t\ \pik, \\
 a_K^{TT}(x) &= \sqrt{K}\int_\bx t^2\ \pik, \\
 \Delta_K^{\varepsilon}(x) &= \sqrt{K}\int_{|(y-x,t)|>\varepsilon}\pik.
\end{align*}

\newcommand{\Ibx}{\mathbf{1}_{B(1)}}

In the calculations that follow, we focus on the first step in our Markov
chain.  We introduce two more random times between $0$ and $\tk_1$, defined
by
\[ S_1=\sup\{t<\tk_1: Y(t)=\xki(t)\}, \]
\[ S_2=\inf\{t>0: Y(t)=\xkii(t)\}. \]
The typical case will be that $0<S_1<S_2<T_1$.  
Lemma \ref{lem:txistinfty} below makes this 
precise.

We also introduce the set $B(1)=\{\omega: |\byk(1)-\byk(0)|<1\}$,
which is the domain of the integrals above.

\newcommand{\tauki}{\tau^K_\infty}
We define $\tauki$ as defined in the second section, to be
the limit of $\lki(t)$ as $t\rightarrow\infty$ in the absence of the
process $\xkii(t)$.  Applying Theorem \ref{thm:escape}, we compute
\[ P(\tauki>t)=\exp(-Kt^2). \]

\newcommand{\taukx}{\lki(S_1)}

\begin{lem}\label{lem:lisismall}
\[ \lim_{K\rightarrow\infty}\sqrt{K}P(\taukx>K^{-1/2+\delta})=0, \]
uniformly in $x$.
\end{lem}
\begin{proof}
This follows from the inequality $\taukx\leq\tauki$ and the explicit
distribution for $\tauki$ (from Theorem \ref{thm:escape}):
\[ \sqrt{K}P(\tauki>K^{-1/2+\delta})
 =\sqrt{K}\exp(-K^{2\delta}). \]
\end{proof}

Next we need to show that $S^K_1$ 
is sufficiently small.  We do this first by examining the duration
of excursions $\xki$ makes from the path of the inert particle.
The measures are from Theorem \ref{thm:exdens}.

\begin{lem}
Define $A^K_\varepsilon$ to be the number of excursions, 
of duration $\varepsilon$ or larger, that $\xki$ makes 
from $\yk$ before time $T_1$.
\[ \lim_{K\rightarrow\infty}\sqrt{K}
  P(A^K_\varepsilon>0 ; \taukx<K^{-1/2+\delta})=0. \]
\end{lem}
\begin{proof}
$A^K_\varepsilon=\#\{(l,t)\in(\varepsilon,\infty)\times[0,\taukx)\}$.
If we condition the process $\xki$ not to make an infinite duration
excursion from $Y$, then $A^K_\varepsilon$ is a Poisson
random variable with rate bounded above by
\[ \int_{\tau=0}^{K^{-1/2+\delta}}\int_{l=\varepsilon}^\infty
  \frac{e^{-K^2\tau^2 l/2}}{\sqrt{2\pi l^3}}\ dl\ d\tau. \]
By a change of variables, we get
\begin{align*}
 &= \int_{l=\varepsilon}^\infty \frac{1}{Kl^2\sqrt{\pi}}
 \int_{u=0}^{K^{1/2+\delta}\sqrt{l/2}}
  e^{-u^2}\ du\ dl \\
 &< \frac{1}{K\sqrt\pi}\int_{l=\varepsilon}^\infty \frac{1}{l^2}
 \int_{u=0}^\infty
  e^{-u^2}\ du\ dl \\
 &= \frac{2}{K\varepsilon}.
\end{align*}
Then 
\[\sqrt{K}P(A^K_\varepsilon>0 ; \taukx<K^{-1/2+\delta})
 \leq \sqrt{K}(1-e^{-2/K\varepsilon}), \]
which converges to $0$ as $K\rightarrow\infty$.  In fact, this will converge
to $0$ for $\varepsilon(K)=K^{-1/2+\delta}$, a fact we will use for the
next lemma.
\end{proof}

\begin{lem}\label{lem:t1small}
\[\sqrt{K}E(S^K_1\Ibx)\rightarrow 0\]
\end{lem}
\begin{proof}
By the previous lemma, we need only consider excursions of length less than
$K^{-1/2+\delta}$ on the set where $\taukx<K^{-1/2+\delta}$.  Then
\begin{align*}
  \sqrt{K}E(S^K_1\Ibx) &\leq
  \sqrt{K}\int_{l=0}^{K^{-1/2+\delta}}\int_{\sigma=0}^{K^{-1/2}+\delta}
  \frac{le^{-K^2\sigma^2l/2}}
  {\sqrt{2\pi l^3}}\ d\sigma\ dl \\
  & \leq
  \sqrt{\frac{K}{2\pi}}\int_{l=0}^{K^{-1/2+\delta}}\int_{\sigma=0}^{K^{-1/2+\delta}}
  l^{-1/2}\ d\sigma\ dl\\
  & = \frac{\sqrt 2}{\sqrt \pi}K^{-\frac 1 4
    +\frac{3\delta}{2}}\rightarrow 0.
\end{align*}
\end{proof}

The next lemma allows us to work with the much nicer density of
$\tauki$ instead of $\taukx$.

\begin{lem}\label{lem:txistinfty}
\[ \lim_{K\rightarrow\infty}\sqrt{K}P(\taukx<\tau^K_\infty)=0 \]
\end{lem}
\begin{proof}
For $\taukx<\tauki$, the inert particle must cross the gap between 
$X^1_t$ and $X^2_t$ before $S^K_1$, the last meeting time of $X^1_t$ and the
inert particle.  Since the particle is in contact with $X^1_t$ at the
instant $S^K_1$, it is sufficient to show that
\[\lim_{K\rightarrow\infty}\sqrt{K}
  P\left(\sup_{s,t<K^{-1/2+\delta}}(|\xki(s)|+|\xkii(t)-x|)>x\right)
   =0. \]
This is equivalent to showing that
\[\lim_{K\rightarrow\infty}\sqrt{K}
  P\left(\sup_{s,t<K^{-1/2+\delta}}(|B_1(s)|+|B_2(t)|+\taukx)>x\right)
   =0. \]
We bound the LHS by
\[\lim_{K\rightarrow\infty}3\sqrt{K}
  P\left(\sup_{t<K^{-1/2+\delta}}|B_1(t)|>x\right)
   \leq 
   \lim_{K\rightarrow\infty}12\sqrt{K}P\left(B_1(K^{-1/2+\delta})>x\right), \]
which is $0$ by a standard computation.
\end{proof}

We will also need a lower bound for $\taukx$, because the time it takes for
the inert particle to cross the gap between the Brownian particles, $S_2-S_1$,
is approximately $x/K\taukx$.

\begin{lem}
\[ \lim_{K\rightarrow\infty}\sqrt{K}P(\taukx<K^{-(3/4+\delta)})=0 \]
\end{lem}
\begin{proof}
By the previous lemma, we need only show this for $\tauki$.
\begin{align*}
  \lim_{K\rightarrow\infty}\sqrt{K}P(\taukx<K^{-(3/4+\delta)})
  &= \lim_{K\rightarrow\infty}\sqrt{K}(1-\exp(-K^{-1/2-2\delta})) \\
  &= \lim_{K\rightarrow\infty} (1+4\delta)\exp(-K^{-1/2-2\delta}) \\
 &= 0.
\end{align*}
\end{proof}

\begin{lem}\label{lem:ltsq}
\[ \lim_{K\rightarrow\infty} \sqrt{K}E(\taukx)^2\rightarrow 0. \]
\end{lem}
\begin{proof}
Because $\taukx<\tauki$, it is enough to compute the expectation
of $(\tauki)^2$:
\[E(\tauki)^2=\int_0^\infty 2s e^{-Ks^2}\ ds=\frac{1}{K}.\]
Multiplying by $\sqrt K$ and taking the limit yields the result.
\end{proof}

\begin{lem}\label{lem:etau}
\[ \lim_{K\rightarrow\infty}\sqrt{K}E\left(\taukx-\frac{\sqrt\pi}{2}\right)=0. \]
\end{lem}
\begin{proof}
By Lemma \ref{lem:txistinfty}, it is enough to compute the expectation
of $\tauki$:
\[E\tauki=\int_0^\infty e^{-Ks^2}\ ds=\frac{\sqrt\pi}{2\sqrt K}.\]
Multiplying by $\sqrt K$, 
and taking the limit yields the result.
\end{proof}


\newcommand{\limcond}{\lim_{K\rightarrow\infty}\sup_{\frac 1 R <|x|<R}}
\begin{lem}\label{lem:bt}
For all $R>0$,
\[ \limcond\left(
b_K^T(x)-\sqrt\pi x\right)=0. \]
\end{lem}
\begin{proof}
Using Lemma \ref{lem:t1small} we disregard the contribution of $S_1^K$ and 
$T_1^K-S_2^K$.

By \cite[p.~196]{ks:91}, we have that a Brownian motion with drift $\mu$,
starting at $x$, has hitting time density at zero
\[ P^\mu(T\in\ dt)=\frac{x}{\sqrt{2\pi t^3}}
 \exp\left(\frac{-(x-\mu t)^2}{2t}\right)\ dt. \] 
A 
table of integrals (e.g. \cite[p.~353]{bigints}) will reveal that
\[ \int_0^\infty t P^\mu(T\in\ dt) = \frac{x}{\mu}. \]
In our case, assuming $X_1(S_1^K)$ is small, that is, the two Brownian
particles remain close to distance $x$ apart, we have
\[ \sqrt{K}E S_2^K = \sqrt{K}\int_{\tau=K^{-(3/4+\delta)}}^{K^{(-1/2+\delta)}}
   \int_{t=0}^{\infty} t P^{K\tau}(S_2^K\in\ dt)P(\tau^K_x\in\ d\tau) \]
\[ \rightarrow x\sqrt{\pi}. \]
The assumption that $X_1(S_1^K)$ is small can be justified by noting that
$X_1(S_1^K)$ will have mean $L_1(S_1)$ and variance $S_1$, and then applying
Lemmas \ref{lem:lisismall} and \ref{lem:t1small}.
\end{proof}

\begin{lem}
For all $R>0$,
\[ \limcond
a_K^{TT}(x)=0. \]
\end{lem}
\begin{proof}
Using the same densities as in the previous lemma, we compute
\[ \sqrt{K} E (S_2^K)^2 = \sqrt{K}\int_{\tau=K^{-(3/4+\delta)}}^{K^{(-1/2+\delta)}}
   \int_{t=0}^{\infty} t^2 P^{K\tau}(S_2^K\in\ dt)P(\taukx\in\ d\tau) \]
\[ \rightarrow 0 .\]
\end{proof}

\begin{lem}
For all $R>0$,
\[ \limcond\left(
b_K^X(x)-\sqrt\pi\right)=0. \]
\end{lem}
\begin{proof}
The change in $\xkii-\xki$ can be expressed as
\[ \xkii-\xki-x=B_2(T_1)-B_1(T_1)+2\taukx. \]
Taking the expectation leaves
\[ E(\xkii-\xki-x)=2E\taukx, \] 
and the result follows from Lemma \ref{lem:etau}.
\end{proof}

\begin{lem}
For all $R>0$,
\[ \limcond\left(
a_K^{XX}(x)-2x\sqrt{\pi}\right)=0. \]
\end{lem}
\begin{proof}
Following the proof of the previous lemma,
\begin{align*}
(\xkii-\xki-x)^2&=B_2(T_1)^2+B_1(T_1)^2+4(\taukx)^2+2B_2(T_1) B_1(T_1)
 \\ &+4\taukx(B_2(T_1)-B_1(T_1)).
 \end{align*}
Taking the expectation 
leaves
\[ E(\xkii-\xki-x)^2=2 E T_1+4E(\taukx)^2+0+0, \] 
and the result follows from Lemma \ref{lem:ltsq} and Lemma \ref{lem:bt}.
\end{proof}

\begin{lem}
For all $R>0$,
\[ \limcond\left(
a_h^{XT}(x)\right)=0. \]
\end{lem}
\begin{proof}
As in the preceding lemmas,
\[ (\xkii-\xki-x)T_1=(B_2(T_1)-B_1(T_1)+2\taukx)T_1. \]
Taking the expectation leaves
\[ E(\xkii-\xki-x)T_1=2E\taukx T_1. \] 
As in the previous lemmas, we discard the contribution from $T_1-S_2$.
Using the probability densities from Lemma \ref{lem:bt}, it is easy to see
that
\[ \int_0^\infty\int_0^\infty
  t\tau P^{K\tau}(S^K_2\in dt)P(\tau^K_x\in d\tau) = \frac{x}{K}, \]
and the result follows.
\end{proof}

We define the process $\hat\byk(t)$ to
be a piecewise linear process derived from the Markov chain
$\byk(j)$ so that $\hat\byk(j/\sqrt{K})=\byk(j)$.

\begin{thm}
The limiting process $\displaystyle\lim_{K\rightarrow\infty}\hat\byk(t)$ 
has generator
\[ L = X\sqrt\pi\frac{\partial^2}{\partial X^2}
        +\sqrt\pi\frac{\partial}{\partial X} 
        +X\sqrt\pi\frac{\partial}{\partial T}.
\]
\end{thm}
\begin{proof}
Apply \cite[Thm. 11.2.3]{sv} with the preceding lemmas.
\end{proof}

From the $X\sqrt\pi(\partial/\partial T)$ term above, 
we can see that the space-time process $\hat\byk(t)$ runs at a different 
rate than do the original Brownian motions which defined our process.
We can perform a change of time to restore the original clock, by dividing
the generator by $X\sqrt\pi$.

\begin{thm}
The limit of the process $(X_2(t/2)-X_1(t/2))$ as $K\rightarrow\infty$ 
is the 2-dimensional Bessel process.
\end{thm}
\begin{proof}
We actually change the clock by the factor $2X\sqrt\pi$ to get the correct
Brownian motion term, because the original process is the difference of two
Brownian motions.
The generator of the space-time process, after the change of time, is
\[ L = \frac 1 2\frac{\partial^2}{\partial X^2}
        +\frac{1}{2X}\frac{\partial}{\partial X}
        +\frac 1 2\frac{\partial}{\partial T}.
\]
Since the process now has a linear clock rate, the first coordinate of the
process will be the original $(X_2(t/2)-X_1(t/2))$, with the generator $L$
with the $T$ term omitted.  This is exactly the generator of the 
two-dimensional Bessel process.
\end{proof}


\newcommand{\Rd}{\Re^d}
\newcommand{\bdyD}{\partial D}
\newcommand{\bd}{B_{\delta}} 
\renewcommand{\vk}[1]{\left<#1\right>}
\newcommand{\td}[1]{\tilde{#1}}
\newcommand{\wtd}[1]{\widetilde{#1}}
\newcommand{\lk}{L}

\section{A Process with Inert Drift in $\Rd$}

We begin this section by recalling some results by P.~Lions and
A.~Sznitman from \cite{lsznit}.  Let $D$ be an open set in $\Re^d$.
Let $n$ be the inward unit vector field on $\bdyD$.  
We make the following assumptions:
\begin{gather}
\label{eq:szi}\exists C_0, \forall x\in\bdyD, \forall x'\in\bar{D}, 
  \forall \lk\in n(x), (x-x',\lk)+C_0|x-x'|^2\geq 0 \\
\begin{split}\label{eq:szii}
\forall x\in\bdyD, &\text{\ if\ } \exists C\geq 0, \exists \lk\in\Re^d:
  \forall x'\in\bar{D}, (x-x',\lk)+C|x-x'|^2\geq 0, \\ 
  & \text{\ then\ }
  \lk=\theta n(x) \text{\ for some\ } \theta\geq 0
\end{split}
\end{gather}
A domain $D$ is called admissible if there is a sequence $\{D_m\}$ of
bounded smooth open sets in $\Re^d$ such that
\begin{enumerate}
\item $D$ and $D_m$ satisfy (\ref{eq:szi}) and (\ref{eq:szii}),
\item if $x_m\in\bar D_m$, $x_m\rightarrow x$, then $x\in \bar D$, and
\item if $K\subset D$ is compact, then $K\subset D_m$ for $m$ large enough.
\end{enumerate}
We will call a pair $(x_t, \lk_t)$ a solution of the Skorohod problem
$(w, D, n)$ if the following hold:
\begin{enumerate}
\item $x_t\in C([0,\infty),\bar D)$,
\item $\lk_t\in C([0,\infty),\Re^d)$, and $\lk_t$ has
bounded variation on every interval $(0,T)$ for all $T<\infty$,
\item $\displaystyle |\lk|_t=\int_0^t 1_{(x_s\in\bdyD)}\,d|\lk|_s$,
\item $\displaystyle \lk_t=\int_0^t n(x_s)\,d|\lk|_s$,
\item\label{it:it} $x_t=w_t+\lk_t$ for $t\geq 0$.
\end{enumerate}


Notationally, $x^j$ or $(x)^j$ will denote the $j$-th component of $x$
when $x$ is a vector or vector function, and $\tilde x$ will denote
the vector $(x^1,\cdots,x^{d-1})$.

We will call a function $x_t$ a solution to the extended Skorohod
problem if condition \ref{it:it} above is replaced by
\begin{equation}\label{eq:esk}
x_t=w_t+\int_0^t \lk_s\ ds+\lk_t.
\end{equation}

\subsection{Existence and uniqueness when $D$ lies above the graph of a function}

The results in this section will be very similar to the
one-dimensional case.  We assume that $D=\{x\in\Rd:
x^d>f(x^1,\cdots,x^{d-1})\}$, with $f(0,\cdots,0)=0$, and that there
is an $0<\alpha<1$ so that $|f(x)|<1-\alpha$ and
$n(x)^d>\alpha$ for all $x$.

\begin{lem}
  If $x_t=w_t+\lk_t$ is a solution to the Skorohod problem in $D$, then
\begin{equation}\label{eq:ineq}
\alpha|\lk|_t<\lk^d_t<|\lk|_t<|\lk|_t/\alpha.
\end{equation}
\end{lem}
\begin{proof}
  From \cite{lsznit} we have that
\begin{align*}
  \lk_t &= \int_0^tn(x_s)\ d|\lk|_s \intertext{so that} \lk^d_t &=
  \int_0^t(n(x_s))^d\ d|\lk|_s > \alpha|\lk|_t.
\end{align*}
Clearly, $\lk^d_t$ is an nondecreasing function.  Since $\lk_0=0$, we also
have that $|\lk_t|<|\lk|_t$.  Combining these, we get \eqref{eq:ineq}.  A
similar computation shows that
\begin{equation*}
|\widetilde{n(x_s)}|\leq\sqrt{1-\alpha^2}|\lk|_t.
\end{equation*}
\end{proof}

\begin{lem}\label{lem:uniq}
  If $x_t$ and $x'_t$ are two solutions to 
the extended Skorohod problem, and $\exists\varepsilon>0$ such that
\begin{equation}\label{eq:assumend}
\sup_{\stackrel{x,y\in\bdyD}{|x-y|=r}} |n(x)-n(y)|
  <r/\varepsilon\qquad\forall r>0,
\end{equation}
then  $x_t=x'_t$ for all $t$.
\end{lem}
\begin{proof}
  Uniqueness seems to require the further assumption \eqref{eq:assumend}
about $D$, similar
  to the assumption about $\mu$ in the one--dimensional case.  Define
\begin{equation*}
\sigma(r)=\sup_{\stackrel{x,y\in\bdyD}{|x-y|=r}} |n(x)-n(y)|.
\end{equation*}
We require that there is some $\varepsilon>0$ such that
$\varepsilon\sigma(r)<r$ for all $r>0$.  This holds for $C^2$ domains.

We may assume that $x_0$ and $x'_0$ start in $\bdyD$, since uniqueness
is clear until the first hitting time of the boundary.  Then we have
that
\begin{align*}
  \wtd{x_t-x'_t} &=
  \wtd{\lk_t-\lk'_t}+\int_0^t(\wtd{\lk_s-\lk'_s})ds+\wtd{x_0-x'_0}
  \\
  &=\int_0^t\td n(x_s)d|\lk|_t-\int_0^t\td
  n(x'_s)d|\lk'|_t+\int_0^t(\wtd{\lk_s-\lk'_s})ds+\wtd{x_0-x'_0}. \\
  \intertext{From the last line, we can compute that}
  \sup_{r,s<t}|\lk_r-\lk_s| &\leq
  2(|\lk|_t+|\lk'|_t)\sigma(\sup_{r,s<t}|\wtd{x_r-x'_s}|)
  \intertext{so that we get the inequality}
  |\wtd{x_t-x'_t}| &\leq
  2(|\lk|_t+|\lk'|_t)\sigma(\sup_{r,s<t}|\wtd{x_r-x'_s}|)(1+t)
  +|\wtd{x_0-x'_0}|.
\end{align*}
We select $T>0$ so that $2(|\lk|_T+|\lk'|_T)(1+T)<\varepsilon/2$.  Let
$\Delta_t=\sup_{r,s<t}|\wtd{x_r-x'_s}|$.  Then we have that
\[ \Delta_T\leq \frac\varepsilon 2\sigma(\Delta_T)+\Delta_0<\Delta_T/2
+\Delta_0, \]
so $\Delta_T\leq 2\Delta_0$.
Then we get that $|\wtd{x_T-x'_T}|\leq 2|\wtd{x_0-x'_0}|$.  If
$x_0=x'_0$, then $x_t=x'_t$ for $0\leq t\leq T$.  The argument is
completed by restarting the process at time $T$.
\end{proof}

\begin{lem}\label{lem:ktub}
  If $x_t=(w_t+I_t)+\lk_t$ solves the Skorohod problem in $D$ for
  ($w_t+I_t$) and $I_t^d$ is an increasing function with $I^d_0=0$,
  then
\begin{equation}\label{eq:skndeq1}
\lk^d_t \leq \sup_{0\leq s\leq t}(-w^d_s\wedge 0)+(1-\alpha)
\end{equation}
\end{lem}
\begin{proof}
The proof is by contradiction.  Let $T$ be a time such that 
(\ref{eq:skndeq1}) does not hold, and let $S$ be the largest time less
than $T$ such that (\ref{eq:skndeq1}) does hold.  Then for $S<t\leq T$,
\[
(w_t+I_t+\lk_t)^d>x^d_t+I^d_t+\sup_{0\leq s\leq t}(-w^d_s\wedge 0)+(1-\alpha)
>1-\alpha\geq f(\tilde x_t).
\]
This implies that $\lk^d_t$ is constant on the interval $[S,T]$, a contradiction.
\end{proof}

\begin{thm}
  Given a continuous $w_t$ with $w_0\in D$, with $D$ satisfying the
  conditions at the beginning of the section, there is a unique $x_t$
  satisfying \eqref{eq:esk}.
\end{thm}
\begin{proof}
  We will construct a solution as in the one-dimensional case.
  
  We combine the results of Lions and Sznitman with Lemma
  \ref{lem:ktub} to construct, for any $\varepsilon>0$,
  $x^\varepsilon_t$, $\lk^\varepsilon_t$, and $I^\varepsilon_t$
  satisfying
\begin{align*}
  x^\varepsilon_t &= w_t+I^\varepsilon_t+\lk^\varepsilon_t,
\end{align*}
where $\lk^\varepsilon_t$ is the local time of $w_t+I^\varepsilon_t$,
and where
\begin{align*}
  I^\varepsilon_t &=
  \lk^\varepsilon_{T^\varepsilon_n}(t-T^\varepsilon_n) ,
  \qquad T^\varepsilon_n<t<T^\varepsilon_{n+1}, \\
  T^\varepsilon_n &= \inf\{t>0: |\lk^\varepsilon|_t=n\varepsilon\}.
\end{align*}
By Lemma \ref{lem:ktub} and \eqref{eq:ineq}, the family
$\{\lk^\varepsilon_t\}_{\varepsilon>0}$ is bounded for each $t$, and by
Lemma \ref{lem:ktl} below, is equicontinuous in $t$.  We can therefore
apply the Ascoli--Arzel\`a Theorem and find a subsequence converging
uniformly on any $[0,T]$ to some $\lk_t$ and $x_t$.  The uniform
convergence gives that $x_t$ is a solution to \eqref{eq:esk}, and
Lemma \ref{lem:uniq} shows that it is the only such solution.
\end{proof}

\subsection{Existence and uniqueness when $D$ is a bounded domain}

\begin{lem}\label{lem:ktl}
  Let $x_t=w_t+\lk_t$ be a solution to the Skorohod problem in $D$, and
  let $\varepsilon>0$.  There is a $\delta>0$ so that
  $|\lk|_T-|\lk|_S<\varepsilon$ whenever $diam(w_{[S,T]})<\delta$.
\end{lem}
\begin{proof}
  Let $x\in\bdyD$, and choose coordinates where $x$ is the origin and
  $n(x)=\vk{0,\cdots,0,1}$.  Because $D$ is $C^1$, we can find a
  $0<\rho<\varepsilon$ so that for $y\in\bdyD$ and $|y-x|<3\rho$, we
  have $n(y)^d>2/3$, and $|y^d|<\rho/3$.  Let $\delta=\rho/3$.
  Suppose that $\rho<|\lk|_T-|\lk|_S<2\rho$.  Then
  $|x_T-x_S|<2\rho+|w_T-w_S|<3\rho$, and
  $x^d_T-x^d_S>w^d_T-w^d_S+\int_S^T n(x_s)d|\lk|_s>2\rho/3-\rho/3$, so
  that $x_T\not\in\bdyD$, a contradiction.
\end{proof}

It is left to show that for a more general domain solutions to
\eqref{eq:esk} exist.  We can do this by piecing together graphs of
functions.

\begin{thm}
  Solutions to \eqref{eq:esk} exist for $D$ an admissible, bounded,
  $C^1$ domain which satisfies locally the extra condition in Lemma
  \ref{lem:uniq}.
\end{thm}
\begin{proof}
  The construction is standard.  Divide $D$ into neighborhoods
  $N_1, \ldots, N_m$ which are nice,
in the sense that, under an appropriate rotation of the standard 
coordinate system, each $N_j\cap\bdyD$ is a section of the graph of a function 
$f_j$
satisfying the conditions at the beginning of the previous section.
  Assume that $w_t$ first encounters
  $N_1$.  Construct the domain which lies above the graph of $f_1$ 
and construct
  $x^{(1)}_t$ satisfying \eqref{eq:esk} on this new domain.  Let
  $T_1=\inf\{t: x^{(1)}_t\not\in N_1\}$.  Repeat the process starting at
  $T_1$ for the function $w_t=w_t+\lk_{T_1}+\lk_{T_1}(t-T_1)$.  Continue
  that construction, so that the limit $x_t$ satisfies \eqref{eq:esk}
  on $[0, \lim\limits_{n\rightarrow\infty} T_n]$.
  
  We wish to show that $\lim\limits_{n\rightarrow\infty}T_n=\infty$.
  If not, say $T_n\rightarrow T$, then by Lemma \ref{lem:ktl} we must
  have that $\lim\limits_{t\rightarrow T}|\lk_t|=\infty$.  Then there is
  some $1\leq j\leq d$, $R,S<T$ so that $0<\lk^j_R<\lk^j_t$ for $R<t<S$,
  and $\lk^j_S\geq \lk^j_R+\mathrm{diam}(\bdyD)+\mathrm{diam}(w[0,T])$.  But this
  contradicts that $x_R, x_S\in D$.

\end{proof}




\bibliographystyle{amsplain}
\bibliography{reflect}

\def\cprime{$'$} \def\cprime{$'$}
\providecommand{\bysame}{\leavevmode\hbox to3em{\hrulefill}\thinspace}
\providecommand{\MR}{\relax\ifhmode\unskip\space\fi MR }
\providecommand{\MRhref}[2]{%
  \href{http://www.ams.org/mathscinet-getitem?mr=#1}{#2}
}
\providecommand{\href}[2]{#2}
\begin{thebibliography}{1}

\bibitem{bkaspi}
Martin Barlow, Krzysztof Burdzy, Haya Kaspi, and Avi Mandelbaum, \emph{Variably
  skewed {B}rownian motion}, Electron. Comm. Probab. \textbf{5} (2000), 57--66
  (electronic). \MR{2001j:60146}

\bibitem{bzwh}
Krzysztof Burdzy and David White, \emph{A {G}aussian oscillator}, Electron.
  Comm. Probab. \textbf{9} (2004), 92--95 (electronic). \MR{MR2108855}

\bibitem{bigints}
I.~S. Gradshteyn and I.~M. Ryzhik, \emph{Table of integrals, series, and
  products}, sixth ed., Academic Press Inc., San Diego, CA, 2000, Translated
  from the Russian, Translation edited and with a preface by Alan Jeffrey and
  Daniel Zwillinger. \MR{MR1773820 (2001c:00002)}

\bibitem{ks:91}
Ioannis Karatzas and Steven~E. Shreve, \emph{Brownian motion and stochastic
  calculus}, second ed., Springer-Verlag, New York, 1991. \MR{92h:60127}

\bibitem{knight:01}
Frank~B. Knight, \emph{On the path of an inert object impinged on one side by a
  {B}rownian particle}, Probab. Theory Related Fields \textbf{121} (2001),
  no.~4, 577--598. \MR{1 872 429}

\bibitem{lsznit}
P.-L. Lions and A.-S. Sznitman, \emph{Stochastic differential equations with
  reflecting boundary conditions}, Comm. Pure Appl. Math. \textbf{37} (1984),
  no.~4, 511--537. \MR{85m:60105}

\bibitem{royden}
H.~L. Royden, \emph{Real analysis}, third ed., Macmillan Publishing Company,
  New York, 1988. \MR{90g:00004}

\bibitem{sv}
Daniel~W. Stroock and S.~R.~Srinivasa Varadhan, \emph{Multidimensional
  diffusion processes}, Grundlehren der Mathematischen Wissenschaften
  [Fundamental Principles of Mathematical Sciences], vol. 233, Springer-Verlag,
  Berlin, 1979. \MR{81f:60108}

\end{thebibliography}


\end{document}